\documentclass{article}

\usepackage[english]{babel}
\usepackage{xcolor}

\usepackage[letterpaper,top=2cm,bottom=2cm,left=3cm,right=3cm,marginparwidth=1.75cm]{geometry}
\usepackage{amsmath}
\usepackage{amsfonts}
\usepackage{graphicx}
\usepackage[colorlinks=true, allcolors=blue]{hyperref}
\usepackage{algorithm}

\begin{document}

  \title{\bf Overcoming global sensitivity limitations: using active subspaces to explore discrepancies between global and local parameter sensitivities}
  \author{Huiyan Zou\\
    Department of Mathematical Sciences, Lafayette College\\
    and \\
    Allison L. Lewis\footnote{Corresponding author: lewisall@lafayette.edu} \\
    Department of Mathematical Sciences, Lafayette College}

\maketitle

\abstract{Global sensitivity metrics are essential tools for assessing parameter importance in complex models, particularly when precise information about parameter values is unavailable. In many cases, such metrics are used to provide parameter rankings that allow for necessary dimension reduction in moderate-to-high dimensional systems. However, globally-derived sensitivity results may obscure localized variability in parameter sensitivities, resulting in misleading conclusions about parameter importance and ensuing consequences for subsequent tasks such as model calibration and surrogate model construction. In this study, we illustrate how discrepancies between globally- and locally-based sensitivity information may arise for an emerging sensitivity metric based on active subspace methodology, as well as for other commonly used sensitivity techniques. In response, we outline a framework that exploits the active subspace to evaluate the stability of parameter sensitivities over the admissible parameter space. This analysis allows one to determine the subregions of the parameter space in which a globally-derived sensitivity metric may be trustworthy.  The proposed framework is illustrated on a collection of simple examples for ease of visualization, as well as the highly-applicable Lotka-Volterra model, for which we demonstrate how these issues may be exacerbated in higher dimensions. Our findings suggest that globally-derived sensitivity information should be treated with caution, and that incorporating analysis on local subregions may improve robustness and accuracy in downstream modeling tasks.
}\\

\noindent%

\section{Introduction}

The purpose of sensitivity analysis is to quantify how variations in model inputs affect model outputs. Such analysis is widely used to guide parameter selection, allow for model simplification, or ensure model identifiability. Global sensitivity analysis measures variations in model output due to parameter perturbations across the entire specified parameter space. Common methods for performing global analysis include screening methods such as Morris elementary effects, variance-decomposition methods such as Sobol' indices, FAST, or eFAST, or correlation and regression methods such as partial correlation coefficients, among others \cite{Saltelli, Smith, Morris, Sobol, Cukier}. In contrast, local sensitivity analysis utilizes linearization at a set of nominal parameter values to measure model behavior more precisely at a particular point of interest. 

When prior information about parameter values is limited, researchers often rely on global sensitivity analysis. However, conclusions derived via such a global approach may not accurately reflect model behavior in all regions of the parameter space, especially in nonlinear and/or high-dimensional models \cite{Brouwer}. The inherent averaging of sensitivity information from a random sampling of points across the space may over-inflate parameter importance due to a small region of out-sized influence, or prioritize parameters that are somewhat sensitive everywhere at the cost of parameters that are generally inert but critically sensitive in certain local regions. For example, there exist safety-critical models such as those used in nuclear engineering for which knowledge of local sensitivities is essential; if such a system is highly sensitive in a specific region, relying solely on global results that may obscure these region-specific sensitivities could yield critical risks \cite{Smith}. Global methods are also highly dependent upon the imposed parameter distributions and ranges, and incapable of differentiating between regions of the space that might produce implausible model outputs due to parameter-driven bifurcations. Averaging of sensitivity results from regions with drastically different model behavior can result in discrepancies with regard to parameter importance \cite{Qian, Broeke, Shin}. As such, an overreliance on global sensitivity measures is likely to lead to poor performance in downstream model analysis tasks. 

Each of the global sensitivity methods mentioned above provides a ranking of parameters in accordance with their overall influence on the model output. However, these methods do not provide a systematic way to evaluate how well globally-derived rankings represent parameter sensitivity in more localized regions. In this investigation, we study how sensitivity information derived from emerging active subspace methodology can provide all of the same information as other commonly-used metrics, while also addressing this limitation. The active subspace---computed via an eigendecomposition of a matrix containing gradient vectors sampled at points from around the parameter space---identifies influential directions in the input space, which are defined to be linear combinations of parameters that contribute most to variation in the model output \cite{blue book, Smith}. By prioritizing linear combinations of inputs that need not align with the standard coordinate space, the active subspace addresses not just basic input prioritization but is also ideally suited for examining how parameters may interact to produce large changes in model response. Connections between active subspace methodology and traditional sensitivity methods (specifically use of the sensitivity Fisher information matrix, or sFIM) have been previously studied \cite{Brouwer}. Here, we introduce a framework that utilizes distances between locally- and globally-defined active subspaces to allow the user to quantify the degree to which parameter sensitivity varies across the admissible space. By identifying regions in which global results are not indicative of local behavior, researchers are given a roadmap as to the conditions under which global metrics may be considered trustworthy, and can identify areas of the parameter space that are worthy of further focus.

In Section \ref{sec:methodology}, we introduce the active subspace, define the corresponding metric used to rank parameters in terms of importance, and introduce our strategy to assess the stability of sensitivity results across the parameter space. Section \ref{sec:results} presents a detailed analysis of several simple two-dimensional models for visualization purposes and the highly-applicable Lotka-Volterra model as an example in higher dimensions. Using the Lotka-Volterra model, we then illustrate the potential downstream impacts of an overreliance on global information when performing subsequent tasks such as model calibration and surrogate model construction. Finally, Section \ref{sec:conclusion} concludes with a summary of the findings and suggestions for future work.

\section{Methodology} \label{sec:methodology}

\subsection{Constructing the Active Subspace}
Active subspace methodology enables dimension reduction in complex models by identifying linear combinations of input parameters that most significantly influence the model output, in addition to orthogonal directions that can be safely neglected in ensuing analysis with minimal loss in accuracy. The following provides a concise summary of active subspace construction; for further details, see \cite{blue book, Smith}.

Let $q = f(\mathbf{x})$ represent the mapping from the parameter space to a scalar-valued model quantity of interest (QoI), such that $f:\mathbb{R}^m \rightarrow \mathbb{R}$. The direction of greatest change in the QoI with respect to the input parameters is captured by the function's gradient, $\nabla_{\mathbf{x}} f(\mathbf{x})$. We introduce the matrix $\mathbf{C}$, 
\begin{eqnarray*}
    \textbf{C} = \int \nabla_{\mathbf{x}} f(\mathbf{x}) \nabla_{\mathbf{x}} f(\mathbf{x})^{T} \rho (\mathbf{x}) \,d\mathbf{x},
\end{eqnarray*} which---defined as the average of the outer product of the function's gradient with itself---quantifies how parameters interact with one another to drive changes in the QoI over the parameter density $\rho(\mathbf{x})$. Due to its symmetry, $\mathbf{C}$ admits an eigendecomposition, $\mathbf{C} = \mathbf{W}\Lambda\mathbf{W}^T,$ where $ \Lambda = \text{diag}(\lambda_1, \ldots, \lambda_m) $ contains the the eigenvalues of $\mathbf{C}$ with $ \lambda_1 \geq \ldots \geq \lambda_m \geq 0 $, and the matrix $\mathbf{W}_{m\times m}$ is an orthonormal matrix containing the normalized eigenvectors of $\mathbf{C}$, $ \mathbf{W}_{m\times m} = [\mathbf{w}_1, \ldots, \mathbf{w}_m] $.

The eigendecomposition of $\mathbf{C}$ provides a natural way to define the active subspace. In particular, we separate the eigenvectors of $\mathbf{C}$ into two distinct sets according to a natural split in the magnitude of the eigenvalues, should one occur: 
\begin{eqnarray*}
    \Lambda = \begin{bmatrix} 
\Lambda_1 & \mathbf{0}\\
\mathbf{0} & \Lambda_2 
\end{bmatrix}, \quad
\mathbf{W} = \begin{bmatrix} \mathbf{W}_1 & \mathbf{W}_2 \end{bmatrix},
\end{eqnarray*}
where $\Lambda_1 = \text{diag}(\lambda_1, \ldots, \lambda_n)$ with $n < m$ contains the $n$ eigenvalues that are largest in magnitude, with $\mathbf{W}_1$ containing the corresponding first $n$ eigenvectors. Similarly, $\Lambda_2 = \text{diag}(\lambda_{n+1},\ldots, \lambda_m)$ contains the smaller eigenvalues corresponding to the eigenvectors in $\mathbf{W}_2 = [\mathbf{w}_{n+1},\ldots,\mathbf{w}_{m}]$.
The active subspace is then defined as the space spanned by $\textbf{W}_1$, which captures the most influential directions in the input space with respect to the quantity of interest. By projecting $\mathbf{x}$ into the subspaces defined by the first $n$ eigenvectors (the active subspace) and the remaining $m-n$ eigenvectors (the inactive subspace), respectively, we can define the active variables to be $\mathbf{y}=\mathbf{W}_1^T \mathbf{x}$ and the inactive variables to be $\mathbf{z}= \mathbf{W}_2^T \mathbf{x}$. Thus, any given $\mathbf{x}$ in the parameter space can be expressed in terms of $\mathbf{y}$ and $\mathbf{z}$ as 
\begin{eqnarray*}
    \mathbf{x} = \mathbf{WW}^T \mathbf{x} = \mathbf{W}_1 \mathbf{W}_1^T \mathbf{x} + \mathbf{W}_2 \mathbf{W}_2^T \mathbf{x} = \mathbf{W}_1 \mathbf{y} + \mathbf{W}_2 \mathbf{z}.
\end{eqnarray*}
This decomposition allows us to analyze $ f(\mathbf{x}) $ in a reduced space by considering its behavior with respect to the components $ \mathbf{y} $ and $ \mathbf{z} $. Since $ \mathbf{W}_2 $ contains directions that contribute very little information about the QoI, the function $ f(\mathbf{x}) $ can be approximated by
\begin{eqnarray*}
    f(\mathbf{x}) = f(\mathbf{W}_1 \mathbf{y} + \mathbf{W}_2 \mathbf{z}) \approx f(\mathbf{W}_1 \mathbf{y}),
\end{eqnarray*} which effectively reduces the dimensionality of our model from $m$ inputs to $n$, where $n<m$. 

In moderate-to-high-dimensional models, directly computing the matrix $\mathbf{C}$ is impractical due to the high-dimensional integration required. In practice, $\mathbf{C}$ is approximated using a Monte Carlo algorithm, as outlined in \cite{blue book} and restated in Algorithm \ref{alg:AS}.

\begin{algorithm}[!htb]
\caption{Random sampling to estimate the active subspace \cite{blue book}}
\label{estimate_C}
\begin{enumerate}
    \item Draw $ M $ samples $ \{\mathbf{x}_j\} $ independently according to the density function $ \rho(\mathbf{x})$.
    \item For each $ \mathbf{x}_j $, compute $ \nabla_{\mathbf{x}} f_j = \nabla_{\mathbf{x}} f(\mathbf{x}_j) $.
    \item Approximate
    \[
    \mathbf{C} \approx \hat{\mathbf{C}} = \frac{1}{M} \sum_{j=1}^{M} (\nabla_{\mathbf{x}} f_j)(\nabla_{\mathbf{x}} f_j)^T.
    \]
    \item Compute the eigendecomposition $ \hat{\mathbf{C}} = \hat{\mathbf{W}} \hat{\Lambda} \hat{\mathbf{W}}^T $.
\end{enumerate}
\label{alg:AS}
\end{algorithm}

For models in which the gradient of the QoI is not available analytically, finite differences may be used to estimate $\nabla_{\mathbf{x}}f(\mathbf{x})$. For a derivation of upper bounds on approximation errors due to Monte Carlo sampling and the use of finite differences, see \cite{blue book}.

Traditionally, the choice of dimension $n$ for the active subspace has been based on locating a natural gap in the eigenvalue spectrum, such that the largest gap in the eigenvalues occurs between $\lambda_n$ and $\lambda_{n+1}$. Such a criterion can be shown to lead to the most accurate estimation of the active subspace when using Algorithm \ref{alg:AS} to estimate the eigendecomposition of \textbf{C} \cite{blue book}. However, other methods for identifying the dimension of the active subspace have been explored, including the use of response surfaces, energy thresholds such as those used in principal component analysis (PCA), or error-based criteria to identify possible reductions \cite{Hal,Jol,Col}. The advantages and disadvantages of each method are discussed in \cite{Col}, where the authors examine the tradeoff between ensuring accuracy in the approximated QoI and balancing the need for computational efficiency.  

\subsection{Sensitivity Analysis Using the Active Subspace}

Traditional global sensitivity analysis methods---including screening methods such as Morris elementary effects \cite{Morris} and variance-based decomposition methods such as Sobol' indices \cite{Sobol}---use sensitivity information from across the parameter space to rank the parameters in order of their influence to the QoI. Parameters that are found to be uninfluential are often fixed at their nominal values, such that the reduction of the problem to be studied is aligned with the original coordinate space. We refer to this as \textit{subset}-based reduction. This is in opposition to \textit{subspace}-based reduction, in which the neglected directions may be linear combinations of the original model parameters.

In general, the use of a subset-based metric for performing model reduction is beneficial because it allows for easier interpretability of the results; a subset of parameters is removed from consideration, and the remaining parameters maintain their physical meaning. However, subspace-based reduction---while not as easily interpretable due to its removal of linear combinations of parameters---is often able to retain more of the total information in the model when using the same dimension as the subset-based analysis, due to its optimization of the information content through the rotation of the parameter space.

Construction of the active subspace allows for both subset- and subspace-based dimension reduction. If subset-based reduction is of interest, parameter rankings can be obtained using \textit{activity scores}, derived from the eigendecomposition of the matrix $\mathbf{C}$. The activity score for parameter $i$, denoted $ \alpha_i $, quantifies the overall contribution of the parameter to the QoI \cite{global sensitivity} and is defined as
\begin{eqnarray}
\alpha_i = \sum_{j=1}^m \lambda_j \mathbf{w}_{i,j}^2, \quad  \text{for} \quad i = 1, ..., m, \label{eqn:actscores}  
\end{eqnarray} 
where $\mathbf{w}_{i,j}$ represents the $(i,j)$th element of the matrix $\mathbf{W}$. This sensitivity metric scales each squared eigenvector component by its corresponding eigenvalue. As such, larger $\alpha_i$ values indicate that the parameter in question is contributing to directions of greater importance, thus exerting greater influence on the model quantity of interest.

Throughout this investigation---in addition to reporting the activity scores derived from the active subspace---we will also include parameter rankings based on Morris elementary effects (mean effect $\mu_i$ and absolute mean effect $\mu_i^*$, along with standard deviation $\sigma_i$) \cite{Morris}, and Sobol' indices (first-order indices $S_i$ and total-effect indices $S_{T_i}$) \cite{Sobol}, in order to compare the performance of our active subspace methodology to more traditional sensitivity methods. For each metric, rankings will be derived on both global and local scales. When deviations between global and local rankings occur, subset-based reduction based on global metrics should be used with caution, as it indicates that we may be removing parameters that appear uninfluential globally but are highly sensitive in certain local regions.

Subspace-based dimension reduction is also possible using information from the active subspace, a significant benefit to using this method over others such as Morris screening or Sobol' indices. Specifically, rather than discarding non-influential parameters as determined by the activity scores, one may choose to discard non-influential \textit{directions} in the input space, defined by the inactive variables $\mathbf{z}$. Such a procedure is often used for applications such as surrogate model construction, in which a lower-dimensional model approximation may be constructed on a subspace of the original parameter space in order to increase computational efficiency for future model evaluations. Surrogate model construction utilizing subspace-based reduction will be illustrated in Section \ref{sec:surrogate}. 

\subsection{Analyzing Stability of Sensitivity Information}

Metrics such as activity scores, Morris elementary effects, and Sobol' indices are designed to summarize information about parameter influence at a global scale, but may do so at the expense of obfuscating local sensitivity variations. This issue is particularly prominent in models that contain nonlinear parameter interactions. In particular, the matrix $\mathbf{C}$ which defines the active subspace can be reframed as the integral of the sensitivity Fisher information matrix (sFIM) over the parameter density.  When parameters interact in a nonlinear manner to produce changes in the model output, such averaging will produce different eigenvectors than would be obtained using a local sFIM, and thus lead to a different approximation of the model output \cite{Brouwer}. Our goal is to develop a tool to quantify the degree to which this issue exists in any given model and highlight the consequences of relying solely on globally-derived information when such nonlinearities are present. Here, we present a metric to determine when global and local parameter sensitivities are out of alignment, based on information provided by the active subspace. This framework can help to identify regions of the parameter space in which global information may be trusted or point to subregions in which further analysis should be conducted. For the remainder of this study, any reference to ``local analysis'' will refer to computation of global sensitivity metrics (e.g., activity scores, Morris elementary effects, or Sobol' indices) on a subregion of the admissible parameter space, as opposed to the traditional definition of local sensitivity analysis, in which derivative-based sensitivity information is measured at a single set of nominal values. 

Our subspace-based criterion assesses the stability of the active subspace by measuring its consistency across different regions of the parameter space. Let $ \mathbf{W} = \begin{bmatrix} \mathbf{W}_1 & \mathbf{W}_2 \end{bmatrix}$ and $ \widetilde{\mathbf{W}} = \begin{bmatrix} \widetilde{\mathbf{W}}_1 & \widetilde{\mathbf{W}}_2 \end{bmatrix} $ represent two matrices containing the eigenvectors of the globally-computed and locally-computed matrices $\mathbf{C}$ and $\widetilde{\mathbf{C}}$, respectively, where the local region in consideration is some subregion of the global parameter space. The distance between the two active subspaces defined by the ranges of $ \mathbf{W}_1 $ and $ \widetilde{\mathbf{W}}_1 $ can be calculated as 
\begin{eqnarray}
\text{dist}(\text{range}(\mathbf{W}_1), \text{range}(\widetilde{\mathbf{W}}_1)) = \left\lVert \mathbf{W}_1 \mathbf{W}_1^T - \widetilde{\mathbf{W}}_1 \widetilde{\mathbf{W}}_1^T \right\rVert,
\label{eqn:subdist}
\end{eqnarray}
which quantifies the distance between the projection matrices of the two subspaces \cite{Stewart}. By analyzing the relationship between subspace distance and parameter values, one can determine subregions of the space in which locally- and globally-computed sensitivity information is in close alignment, or identify certain parameter ranges for which the globally-computed active subspace prioritizes different combinations of parameters than would be favored locally. In Section \ref{sec:results}, we will illustrate how such regions can be visually identified for both two-dimensional and higher-dimensional models.  

\section{Results} \label{sec:results} 

In what follows, we illustrate the proposed methodology on several two-dimensional examples---chosen in order to facilitate straightforward visualization of how local and global sensitivities may not align---and on the six-dimensional Lotka-Volterra model, where we demonstrate how such issues may be further intensified in higher dimensions. All of the analysis is performed using MATLAB software.

\subsection{Two-Dimensional Models}

To illustrate the variation in sensitivity that can occur as one moves around the parameter space, we begin with three simple two-dimensional examples. Using $\mathbf{x}\in [0,1]\times[0,1]$, we define our three models as follows:
\begin{eqnarray*}
    f_1(\mathbf{x}) = e^{0.7x_1+0.3x_2}, \ \ \
    f_2(\mathbf{x}) = x_1e^{0.7x_1+0.3x_2},\ \ \ \text{and} \ \
    f_3(\mathbf{x}) = e^{0.7x_1^2x_2+0.3x_2}.
\end{eqnarray*}

\begin{figure}[!htb]
        \includegraphics[width=.32\textwidth]{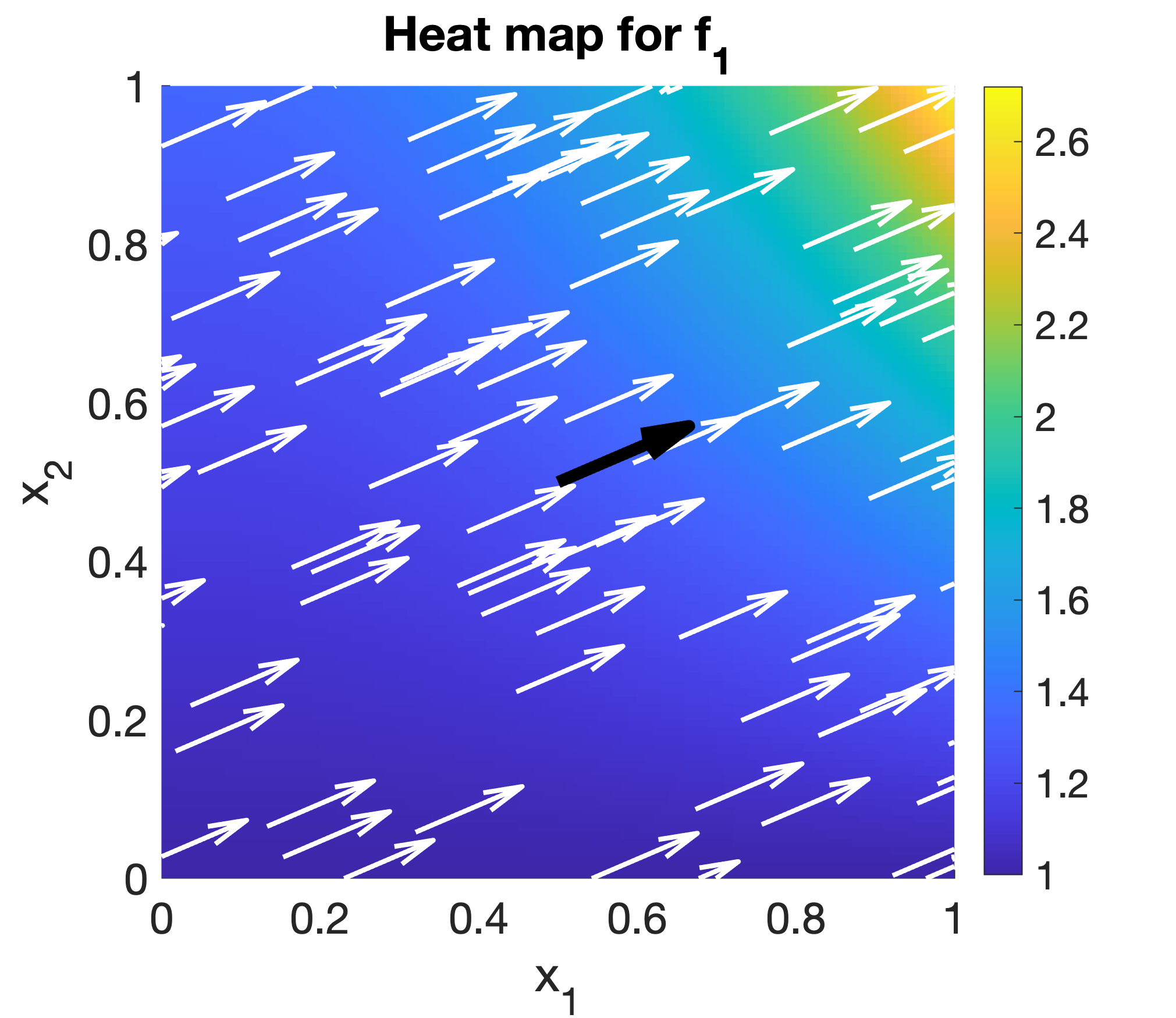}
        \includegraphics[width=.32\textwidth]{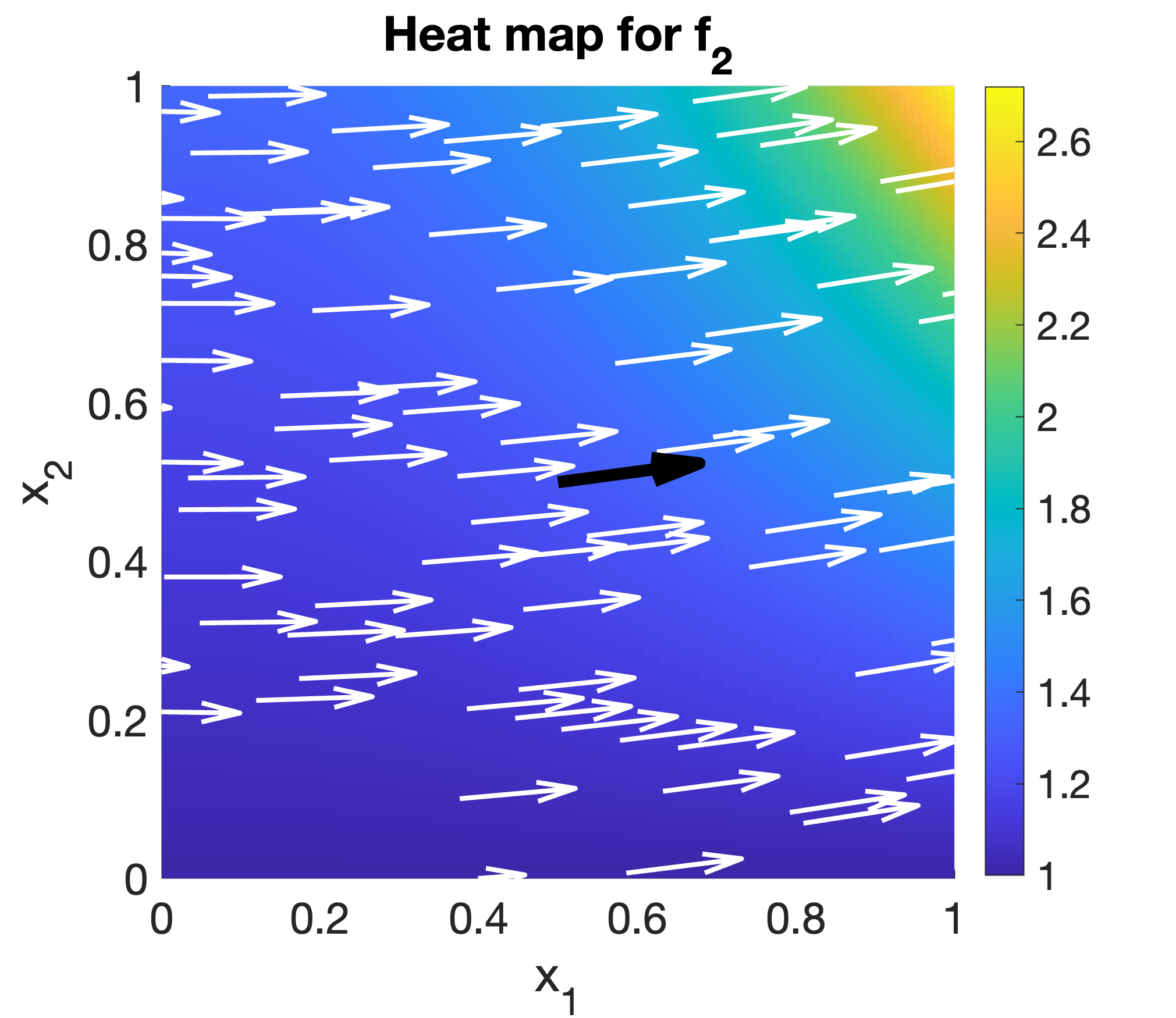}
        \includegraphics[width=.32\textwidth]{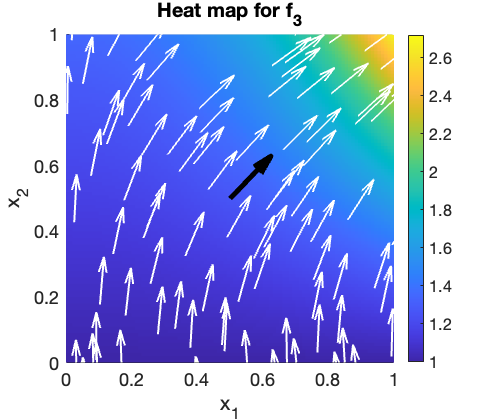}
    \caption{Visualization of locally-computed gradient vectors (white) versus a gradient vector computed from global sampling (black). }
    \label{f1f2heat}
\end{figure}

All three of these examples permit a one-dimensional active subspace---as determined by the negligibility of the second eigenvalue when compared to the first---which is defined by the span of the gradient vector. The variability in gradient vector direction for these examples is visualized in Figure \ref{f1f2heat}. Function $f_1$ exhibits no variation at all over the parameter space; the gradient vector is a scalar multiple of $[0.7,\ 0.3]$ everywhere in the space. For such an example, an active subspace computed globally will be identical to an active subspace computed on any sub-portion of the parameter space, and a ranking of the parameters via activity scores will always prioritize $x_1$ over $x_2$. With function $f_2$, there is a slight amount of variability in the direction of the gradient vector as one moves around the space, though $x_1$ is clearly always more influential than $x_2$. The function $f_3$ is an example of a function for which the direction of the gradient vector is highly dependent on the location in the parameter space where it is computed.

We want to understand how this variation in gradient vector direction affects the subsequent sensitivity analysis. For each function, we divide the domain into local grids of size $0.01 \times 0.01$, and compute the active subspace within each grid square using Algorithm \ref{alg:AS} with $M=10$ samples randomly drawn from a uniform density function using a Latin Hypercube design. The locally-computed active subspaces are then used to determine activity score rankings for each grid square using Equation (\ref{eqn:actscores}). Finally, we use 100,000 samples from a uniform Latin Hypercube design across the full parameter space---corresponding to 10 samples per local region---to compute a global active subspace and global activity score ranking.  

As expected, locally-computed activity scores for functions $f_1$ and $f_2$ indicate that $x_1$ is more influential than $x_2$ in every local region, matching their global rankings. For function $f_3$, the global normalized activity scores of $[0.9412, 1]$ indicate a slight prioritization of $x_2$; however, in 2187 of the 10,000 local regions---namely, those concentrated in the upper right corner---$x_1$ is prioritized over $x_2$ instead. 

\begin{figure}[!htb]
        \includegraphics[width=.32\textwidth]{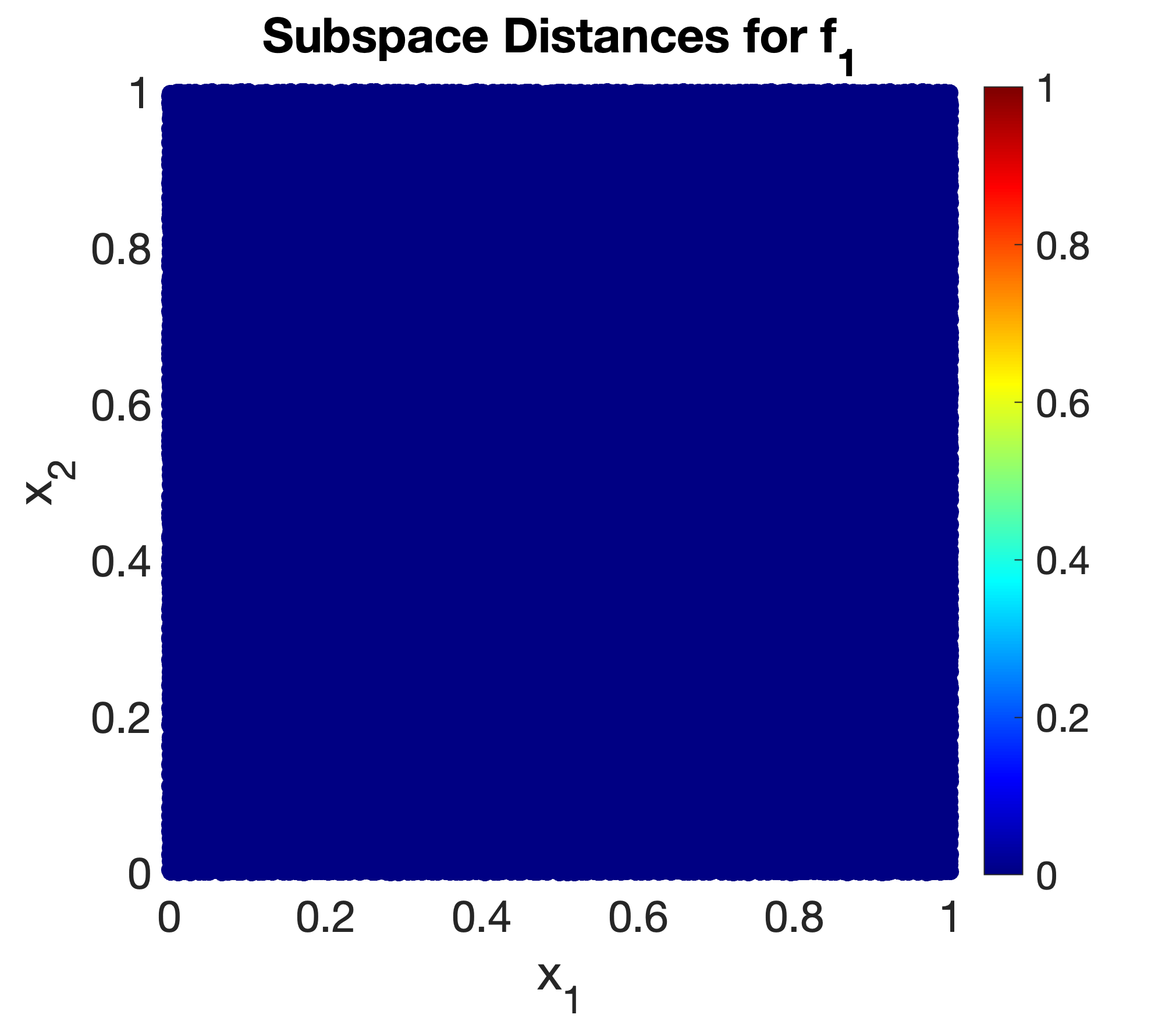}
        \includegraphics[width=.32\textwidth]{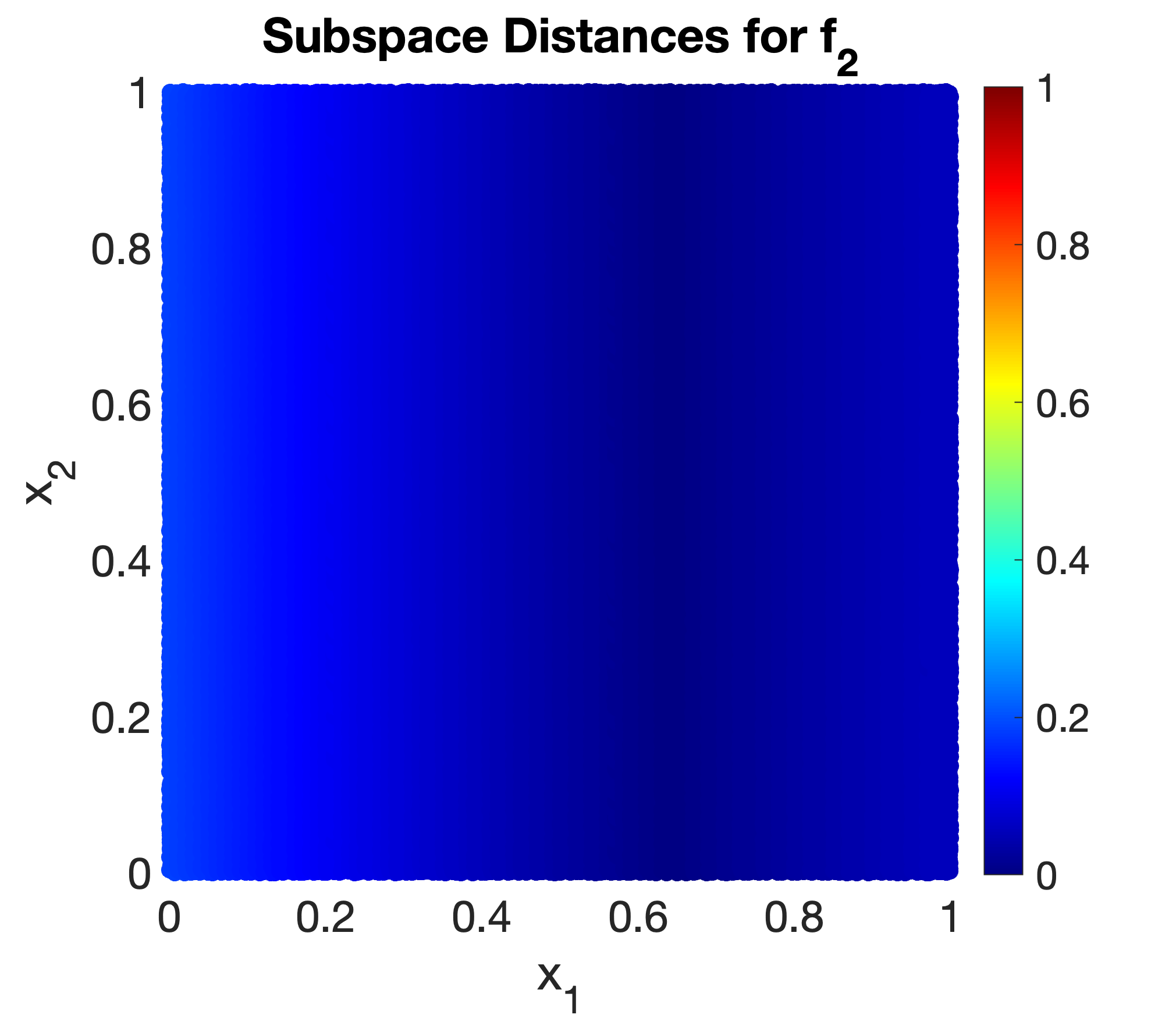}
        \includegraphics[width=.32\textwidth]{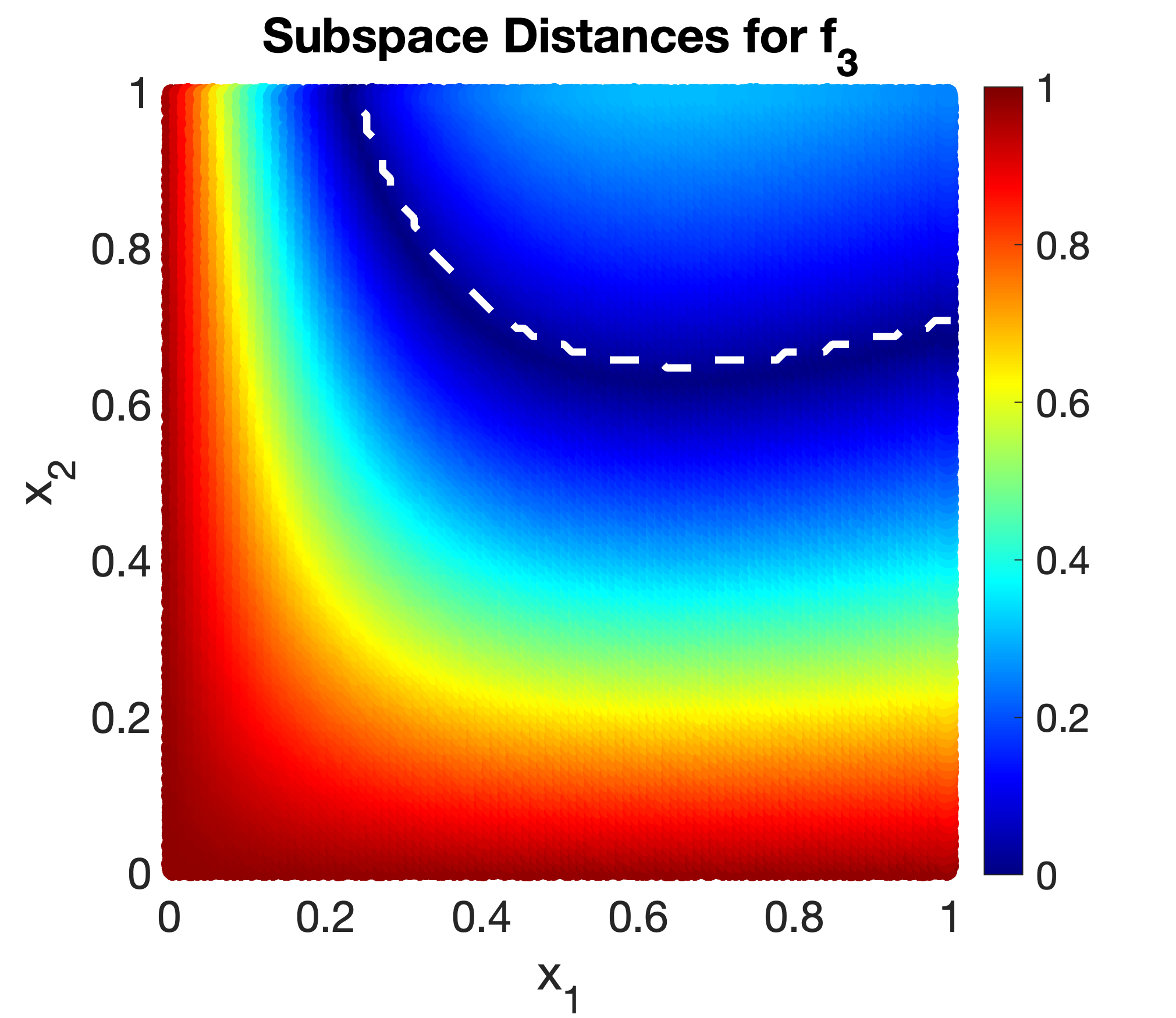}
    \caption{Subspace distance visualization. For function $f_3$, the dashed curve indicates the separation of regions where $x_1$ is ranked first (above) and where $x_2$ is prioritized (below), according to locally-computed activity scores.}
    \label{f1f2}
\end{figure}

Figure \ref{f1f2} illustrates the distances between locally- and globally-computed active subspaces, computed using Equation (\ref{eqn:subdist}). Again, heat maps for $f_1$ and $f_2$ show little or no variation across the space, indicating that the global active subspace is relatively well-aligned with local information everywhere. However, the heat map for $f_3$ indicates that there are regions in the parameter space where the global active subspace is a very poor representation of local activity, most notably when $x_1$ and/or $x_2$ is small: regions where the influence of the parameter $x_2$ strongly outweighs $x_1$. Given these large disparities in active subspace alignment, together with the differing activity score rankings from two regions of the space, this constitutes our first example in which a globally-computed sensitivity metric may give misleading information about function behavior in a local region.

\subsection{Lotka-Volterra Model} \label{LV_section}

\begin{figure}[!htb]
\centering
\includegraphics[width = \textwidth]{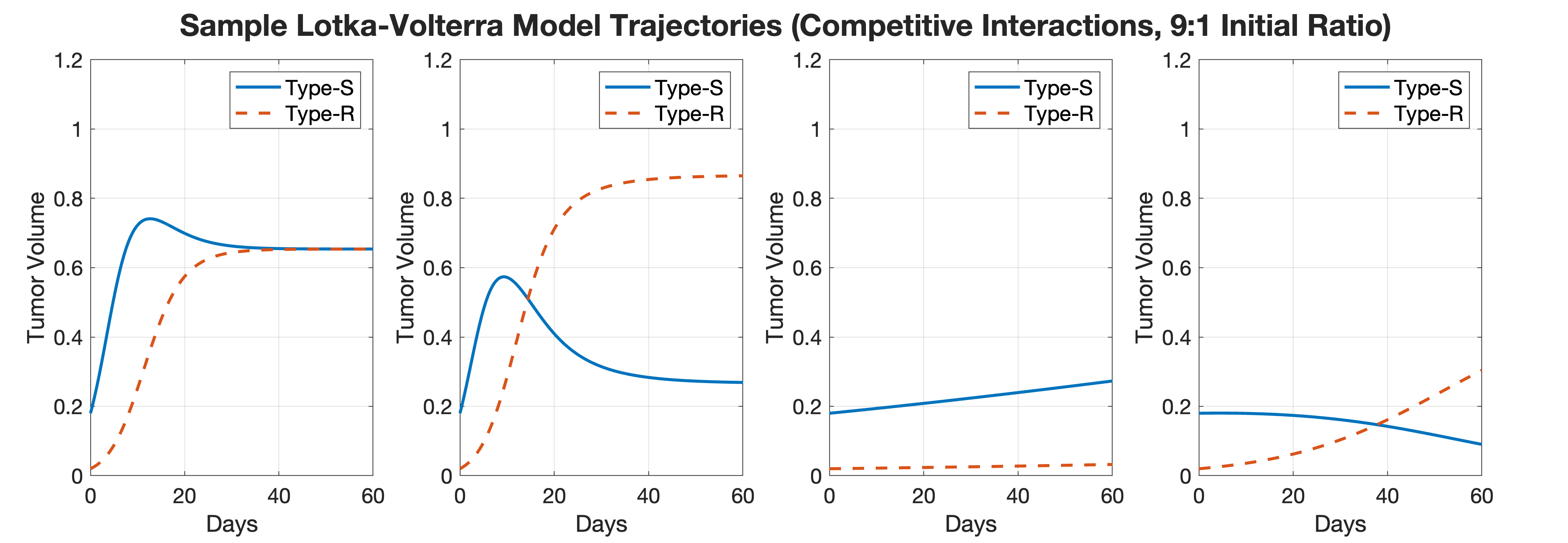}
\caption{Competitive Lotka–Volterra model trajectories with a 9:1 initial Type-$S$ : Type-$R$ ratio. Trajectories for Type-$S$ (blue, solid) and Type-$R$ (orange, dashed) cell volumes illustrate how different choices of growth rates, carrying capacities, and interaction strengths can impact the behavior of the system. }
\label{fig:LV_traj}
\end{figure}

To demonstrate how the disparity in globally- versus locally-derived sensitivity information may be exacerbated in models that are high-dimensional and/or highly nonlinear, we turn our attention to the six-dimensional Lotka-Volterra model. The Lotka-Volterra model is a system of differential equations originally developed to describe predator–prey interactions in ecological systems \cite{Lotka1925,Volterra1926}. It has also been widely adapted across disciplines, including physics \cite{PhysRevB2023,ThreeLevelLaser2021,MDPI2023} and chemistry \cite{ChemPaper1,ChemPaper2}, to describe dynamics in systems such as lasers, reaction kinetics, and phase transitions. The general form of the model captures nonlinear interactions between two or more species or components, making it a flexible structure for exploring the effects of parameter variation in dynamical systems. In this study, we apply our analysis to a variant of the Lotka-Volterra model describing the growth of a heterogeneous tumor spheroid comprised of two cancer cell types, Type-$S$ and Type-$R$ \cite{LV_model}:
\begin{eqnarray*}
\frac{dS}{dt} &=& r_S S\left(1 - \frac{S}{K_S} - \frac{\gamma_R R}{K_S}\right)\\ 
\frac{dR}{dt} &=& r_R R\left(1 - \frac{R}{K_R} - \frac{\gamma_S S}{K_R}\right),
\end{eqnarray*}
where each population is considered to undergo logistic growth in the absence of the other cell line, with growth rates of $r_S$ and $r_R$ and carrying capacities of $K_S$ and $K_R$, respectively. The signs and magnitudes of the parameters $\gamma_S$ and $\gamma_R$ determine interaction strength between the two populations, and $S(t)$ and $R(t)$ represent the volume (mm$^3$) of Type-$S$ and Type-$R$ cancer cells at time $t$.  For this investigation, we initialize the system using a total tumor volume of 0.02 mm$^3$ with a 9:1 initial ratio of Type-$S$ to Type-$R$ cells.  Additionally, we set parameter ranges of $[0,1]$ for all parameters $[r_S,r_R,K_S,K_R,\gamma_S,\gamma_R]$, which restricts us to positive growth scenarios with competitive interactions only. Figure \ref{fig:LV_traj} showcases some of the variability in behavioral patterns that can occur for this system, using different sets of parameter values drawn from within these ranges. We define the scalar-valued QoI to be the integral of the total tumor volume over an eight-week (56-day) growth observance period, 
\begin{eqnarray*}
\int_0^{56} \left(S(t)+R(t)\right) \ dt,
\end{eqnarray*}
 which we approximate computationally using the trapezoidal rule with bin sizes of one day. This QoI can be considered a representation of the total tumor burden for a patient.

For the localized sensitivity analysis, we divide each parameter's range into eight segments of length 0.125, resulting in a total of $8^6 =262,144$ local regions for consideration. Within each region, we generate 10 sample points from a uniform Latin Hypercube design to compute the local active subspace. Global analysis is conducted using $10\times 8^6$ sample points from a Latin Hypercube design, corresponding to 10 sample points per local region. Gradient vectors required for active subspace computations are computed using central finite differences with step size $h= 10^{-5}$, after scaling all parameters to the same range.

Beginning with subset-based analysis, we first compute the globally-determined sensitivity metrics for each parameter in the model, presented in Table \ref{LV_globalmetrics}.  Using the normalized activity scores derived from the active subspace analysis, rankings indicate that the Type-$S$ growth rate, $r_S$, and carrying capacity, $K_S$, have the most significant influence on the QoI (with these two outranking their $r_R$ and $K_R$ counterparts due to the initial cell ratio of 9:1), while the interaction parameters $\gamma_S$ and $\gamma_R$ contribute the least. The overall global ranking is given by $[r_S,\ K_S, \ K_R, \ r_R, \ \gamma_S, \ \gamma_R]$. In contrast, both Morris screening and Sobol' indices prioritize the two carrying capacities, $K_S$ and $K_R$, when analyzed across the full parameter space, with global ranking $[K_S, \ K_R, \ r_R, \ r_S, \ \gamma_S, \ \gamma_R]$. 

\begin{table}[h]
\centering
\begin{tabular}{l|c|cccccc}
    \hline
    \multicolumn{2}{c}{} & $r_S$ & $r_R$ & $K_S$ & $K_R$ & $\gamma_S$ & $\gamma_R$ \\
    \hline
    Activity Scores & $\alpha_i/\alpha_{\text{max}}$ &  1.0000 & 0.6584 & 0.9142 & 0.6777 & 0.0729 & 0.0553 \\
    \hline
    Morris Screening & $\mu$     & 1.0873 & 1.1606 & 4.4601 & 3.9395 & -0.7855 & -0.6298\\
    & $\mu^*$   & 1.1046 & 1.1731 & 4.5060 & 3.9654 & 0.8030 & 0.6442\\
    & $\sigma$  & 2.4914 & 2.3383 & 7.7021 & 8.3844 & 1.0925 & 0.9733\\
    \hline 
    Sobol' Indices & $S_i$     & 0.0834 & 0.0908 & 0.3701 & 0.2435 & 0.0122 & 0.0083 \\
    & $S_{Ti}$  & 0.1513 & 0.1747 & 0.4719 & 0.3472 & 0.0316 & 0.0237 \\
    \hline
\end{tabular}
\caption{Global sensitivity metrics for all methods.}
\label{LV_globalmetrics}
\end{table}

To understand whether global results can be accurately applied to analyze model behavior on a more localized scale, we conduct our analyses individually on the set of $8^6$ locally-defined subregions and look for patterns in which parameters are prioritized. For the active subspace analysis, we observe 306 unique parameter rankings (out of 720 possible rank orders) among the $8^6$ regions, revealing that sensitivity information varies drastically around the space. The Morris screening and Sobol' indices methods reveal the observance of 305 and 720 unique rankings, respectively, illustrating that this issue is not unique to the activity score metric. 

The top ten most frequently observed local rankings for the activity score metric are tallied in Table \ref{tab:LV local ranking}. Notably, the global activity score ranking of $[r_S,\ K_S, \ K_R, \ r_R, \ \gamma_S, \ \gamma_R]$ is only the 53rd most frequently observed local ranking, occurring in only 1526 out of the total 262,144 regions. The two most globally-favored parameters, $r_S$ and $K_S$, are ranked as the top two parameters in only three of the ten most commonly observed local rankings. While $K_S$ is regularly observed among the three most influential parameters in the most frequently observed rankings, $r_S$ is actually ranked as the \textit{least} influential parameter in five of the top ten observed local rankings, casting doubt upon the reliability of the global active subspace findings.

\begin{table}[h]
    \centering
    \begin{tabular}{ccc|ccc}
        \hline
        Order & Ranking & Frequency & Order & Ranking & Frequency \\
        \hline
        1 & $K_R$, $r_R$, $K_S$, $\gamma_S$, $\gamma_R$, $r_S$ & 12,263 & 6 & $K_R$, $K_S$, $\gamma_R$, $r_R$, $\gamma_S$, $r_S$ & 5,734 \\
        2 & $K_S$, $r_S$, $K_R$, $\gamma_S$, $\gamma_R$, $r_R$ & 10,497 & 7 & $K_S$, $K_R$, $\gamma_S$, $r_S$, $\gamma_R$, $r_R$ & 4,923 \\
        3 & $K_S$, $r_S$, $K_R$, $\gamma_S$, $r_R$, $\gamma_R$ & 9,392 & 8 & $K_S$, $K_R$, $\gamma_R$, $\gamma_S$, $r_R$, $r_S$ & 4,758 \\
        4 & $K_S$, $r_S$, $K_R$, $\gamma_R$, $\gamma_S$, $r_R$ & 9,052 & 9 & $r_R$, $K_R$, $K_S$, $\gamma_S$, $\gamma_R$, $r_S$ & 4,629 \\
        5 & $K_R$, $K_S$, $r_R$, $\gamma_R$, $\gamma_S$, $r_S$ & 5,782 & 10 & $K_R$, $r_R$, $K_S$, $\gamma_S$, $r_S$, $\gamma_R$ & 4,533 \\
        \hline
    \end{tabular}
    \caption{Ten most frequently observed local activity score rankings for the Lotka-Volterra model.}
    \label{tab:LV local ranking}
\end{table}

\begin{figure}[!htb]
\centering
\includegraphics[width = \textwidth]{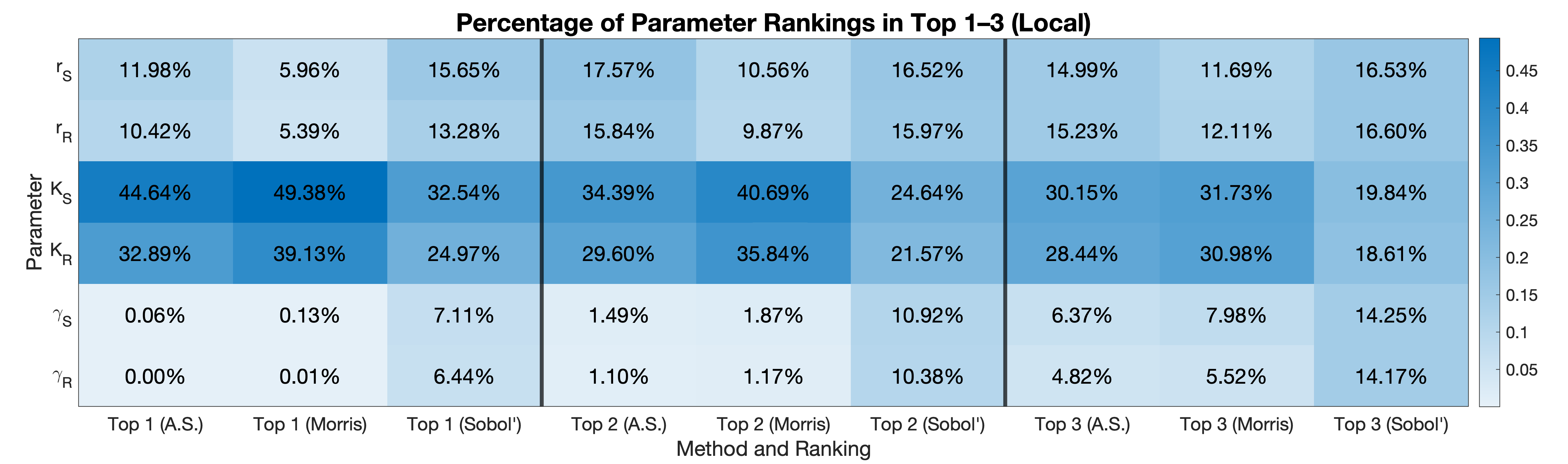}
\caption{Percentage of local regions in which each parameter is ranked among the top 1, top 2, or top 3 of the ranking list.}
\label{Top123}
\end{figure}

Figure \ref{Top123} displays the percentage of cases in which each parameter appears in the top 1, top 2, or top 3 parameter ranking spots among all of the localized regions, for each of the activity score, Morris screening, and Sobol' indices methods. Despite the differing global rankings, all methods are in agreement about which parameters are most often flagged as influential among the locally-defined regions.  For all three methods, one of the two carrying capacities is most often chosen as the most sensitive parameter. The growth rates are chosen as most influential far less frequently but often show up in the top 2 or top 3 of the ranking list. Finally, the interaction parameters are rarely (if ever) chosen as the most influential, but do occasionally make appearances in the top 2 or top 3 rankings---especially in the case of Sobol' indices, where parameter ranking distributions are more evenly dispersed.

\begin{figure}[!htb]
\centering
\includegraphics[width=.5\textwidth]{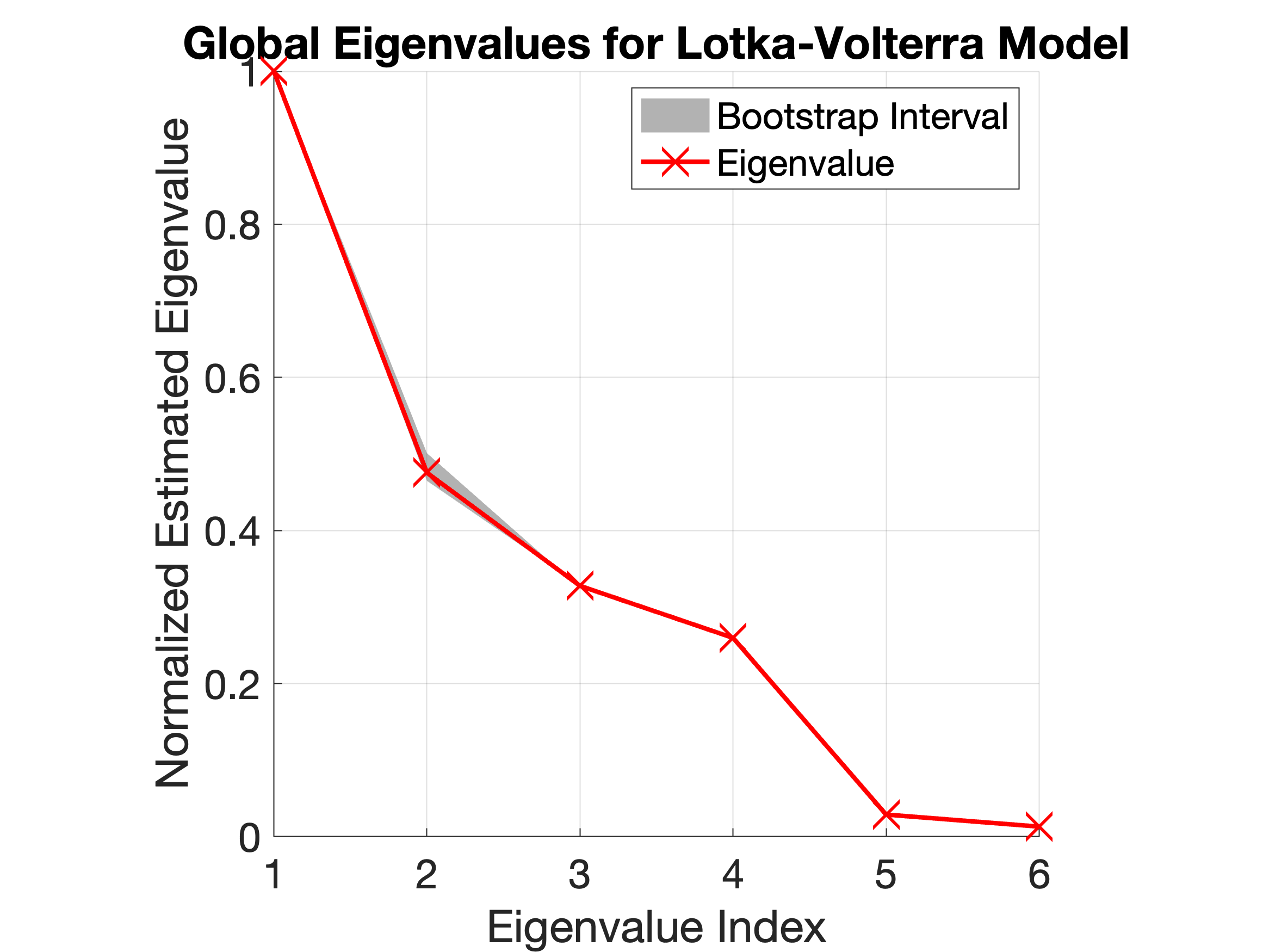}
\caption{Eigenvalues of the globally-computed matrix $\hat{\mathbf{C}}$ for the Lotka-Volterra model.}
\label{LV_eigenval}
\end{figure}

Our subset-based analysis reveals that many different rankings may appear in localized regions that differ from the global results and that this issue persists across multiple established sensitivity methods. To determine in which local regions the global results can be trusted, we conduct a subspace-based analysis to analyze the stability of the active subspace as we move around the admissible parameter space. An analysis of the eigendecomposition of the globally-computed matrix $\mathbf{C}$ suggests the use of a 4-dimensional active subspace; eigenvalue decay is relatively gradual, as seen in Figure \ref{LV_eigenval}, but the final two eigenvalues appear negligible in magnitude compared to the first four. Given this structure, we select the first four eigenvectors of the globally- and locally-computed eigendecompositions to define the matrices $\mathbf{W_1}$ and $ \widetilde{\mathbf{W_1}}$, which are used in Equation (\ref{eqn:subdist}).

\begin{figure}[!htb]
        \includegraphics[width=.32\textwidth]{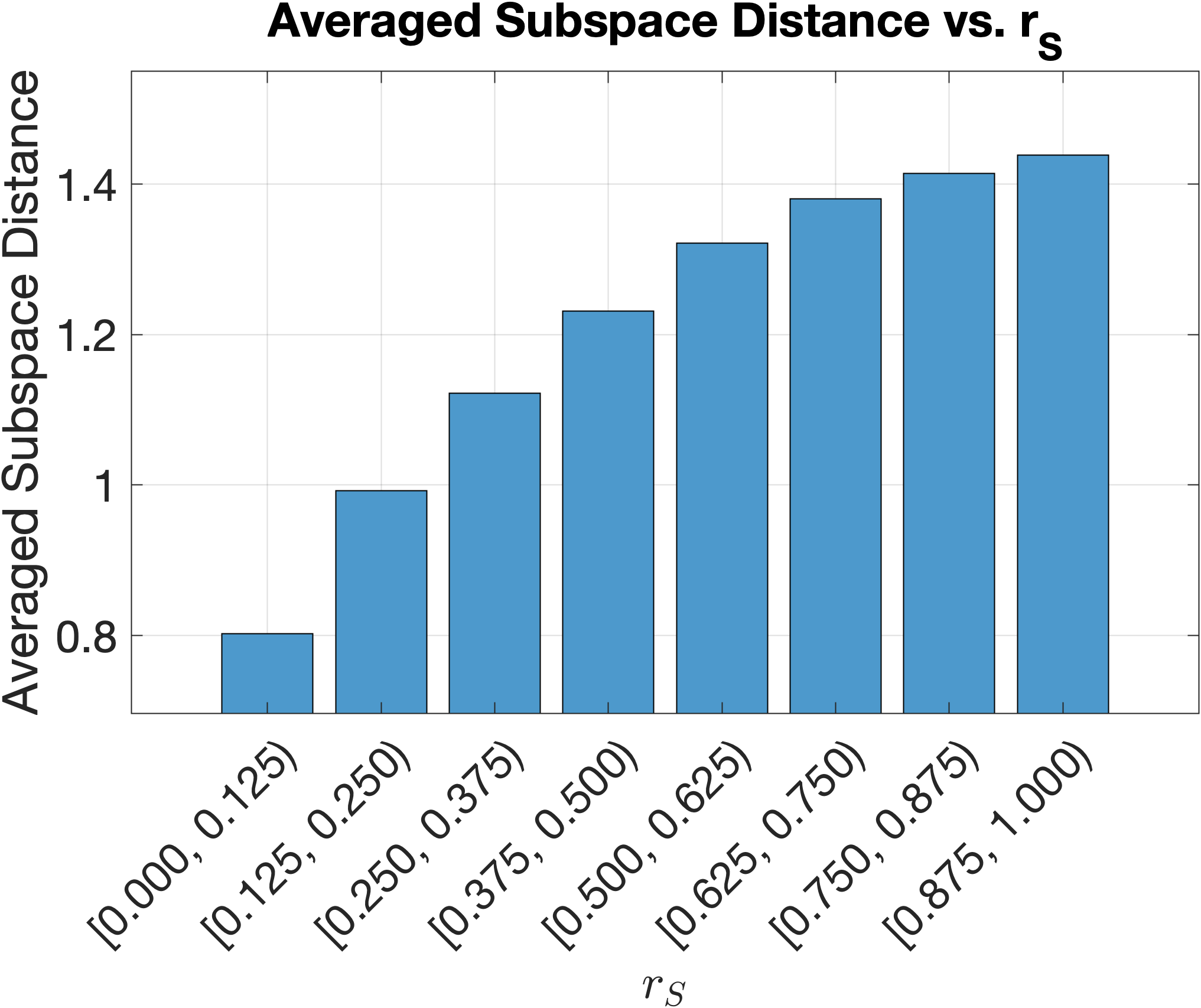}
        \includegraphics[width=.32\textwidth]{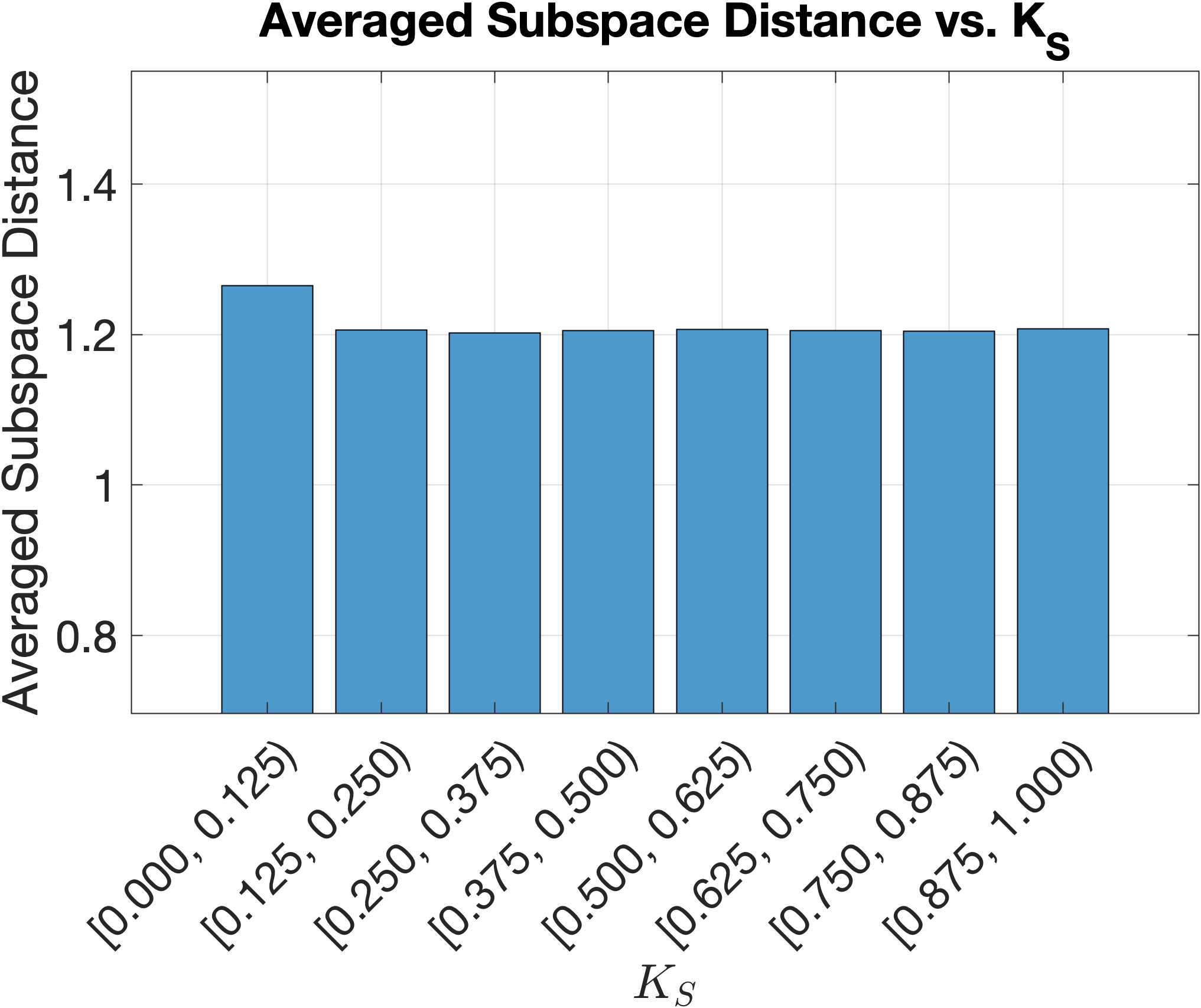}
        \includegraphics[width=.32\textwidth]{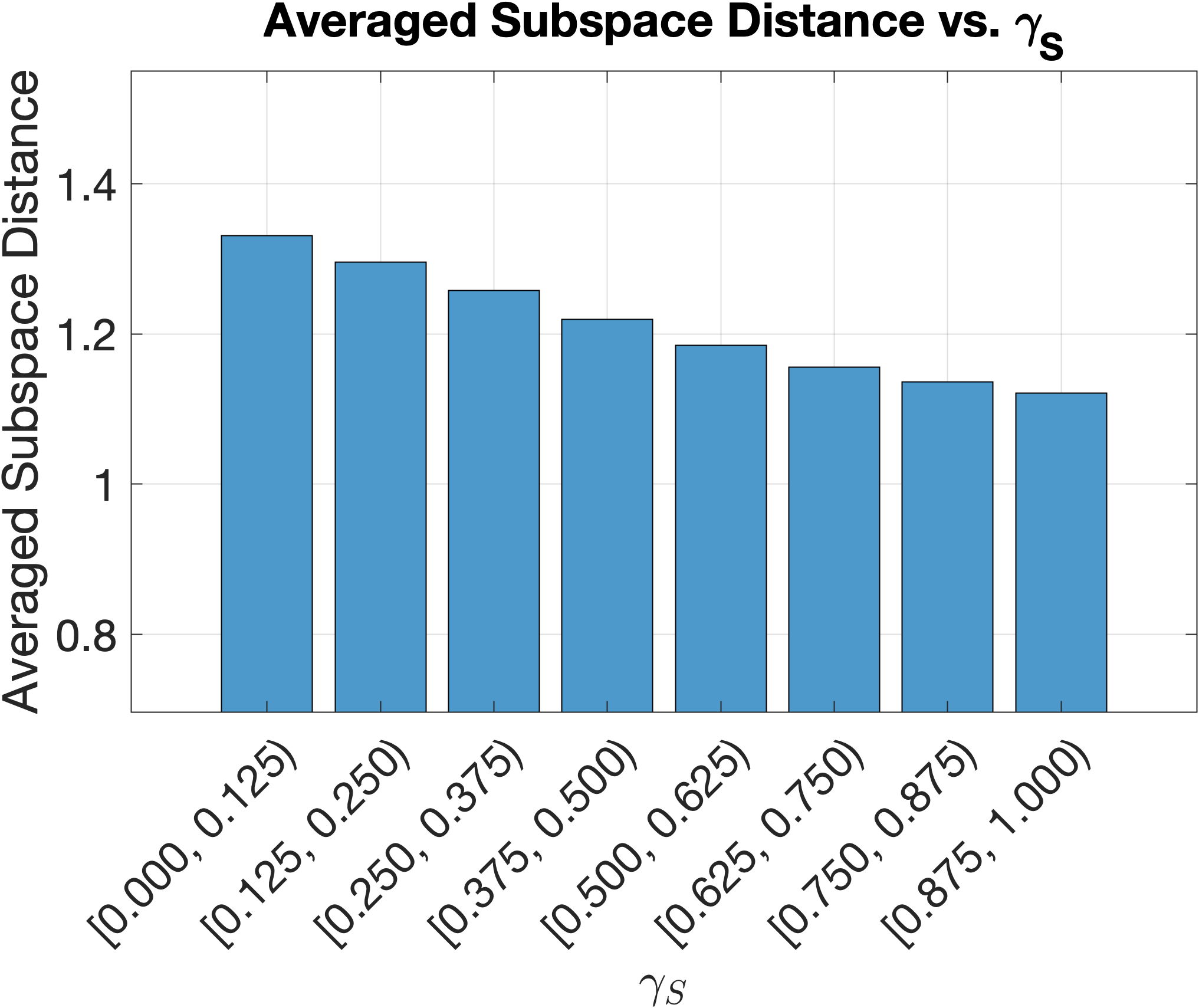} \\ \ \\
        \includegraphics[width=.32\textwidth]{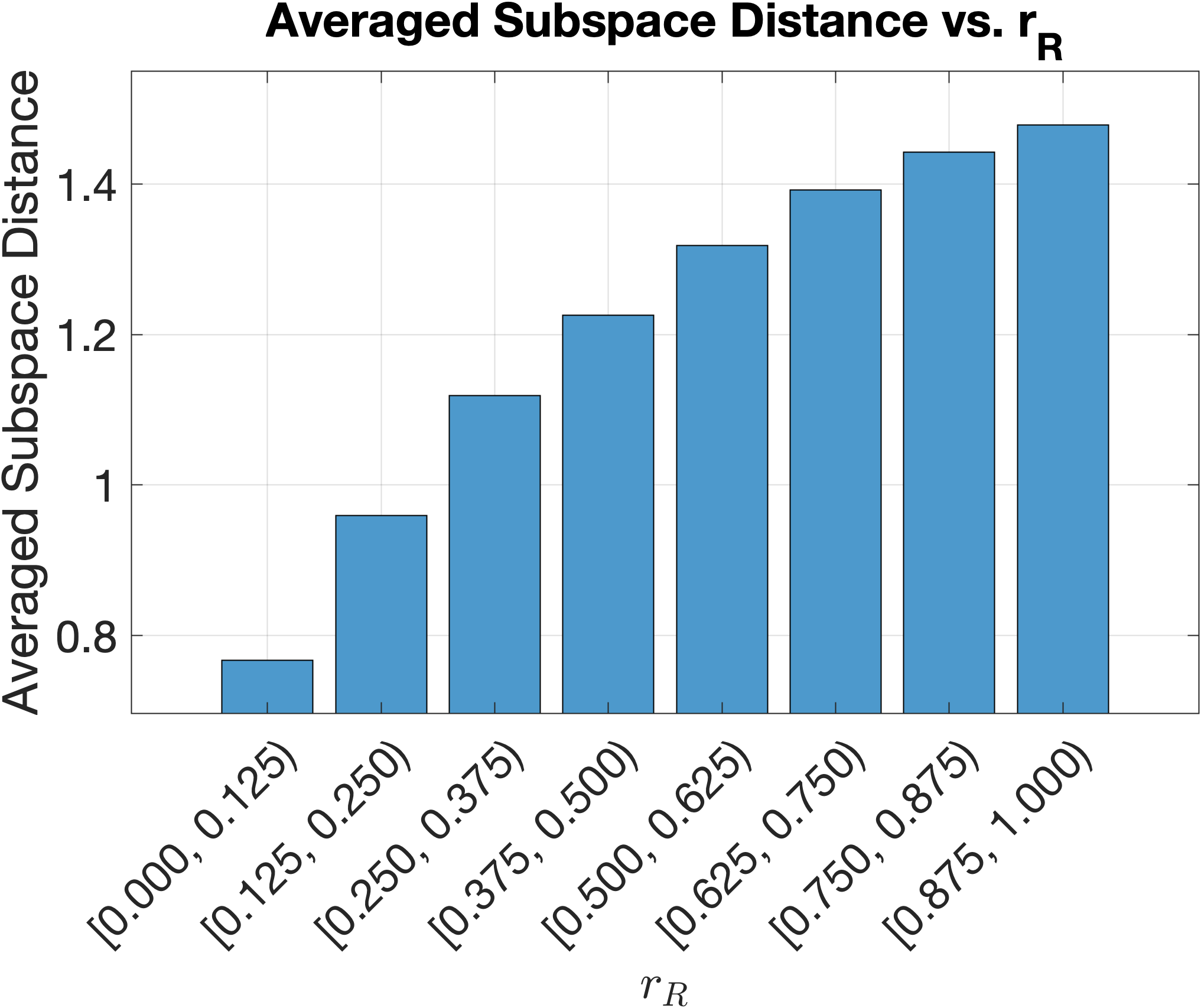}
        \includegraphics[width=.32\textwidth]{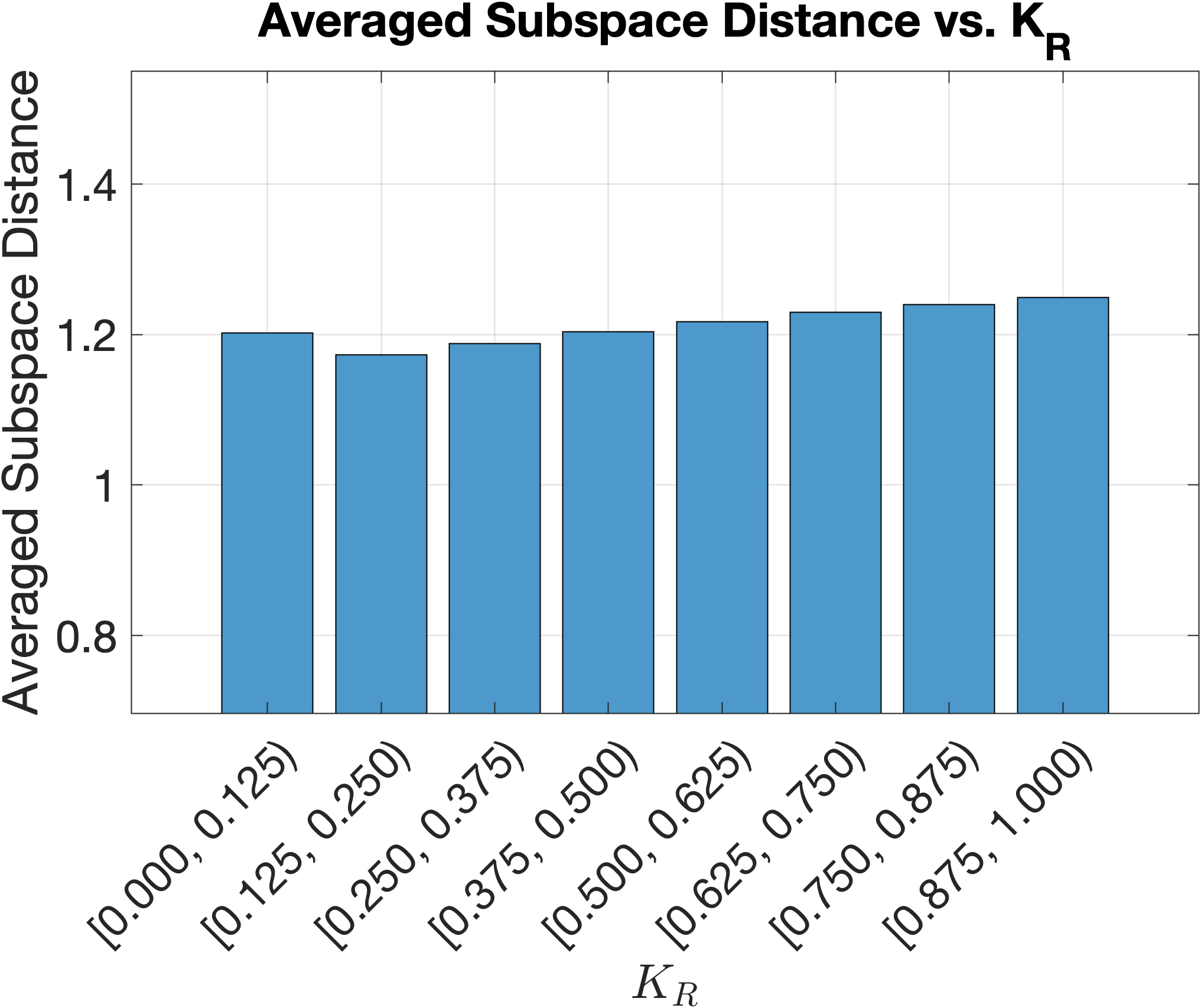}
        \includegraphics[width=.32\textwidth]{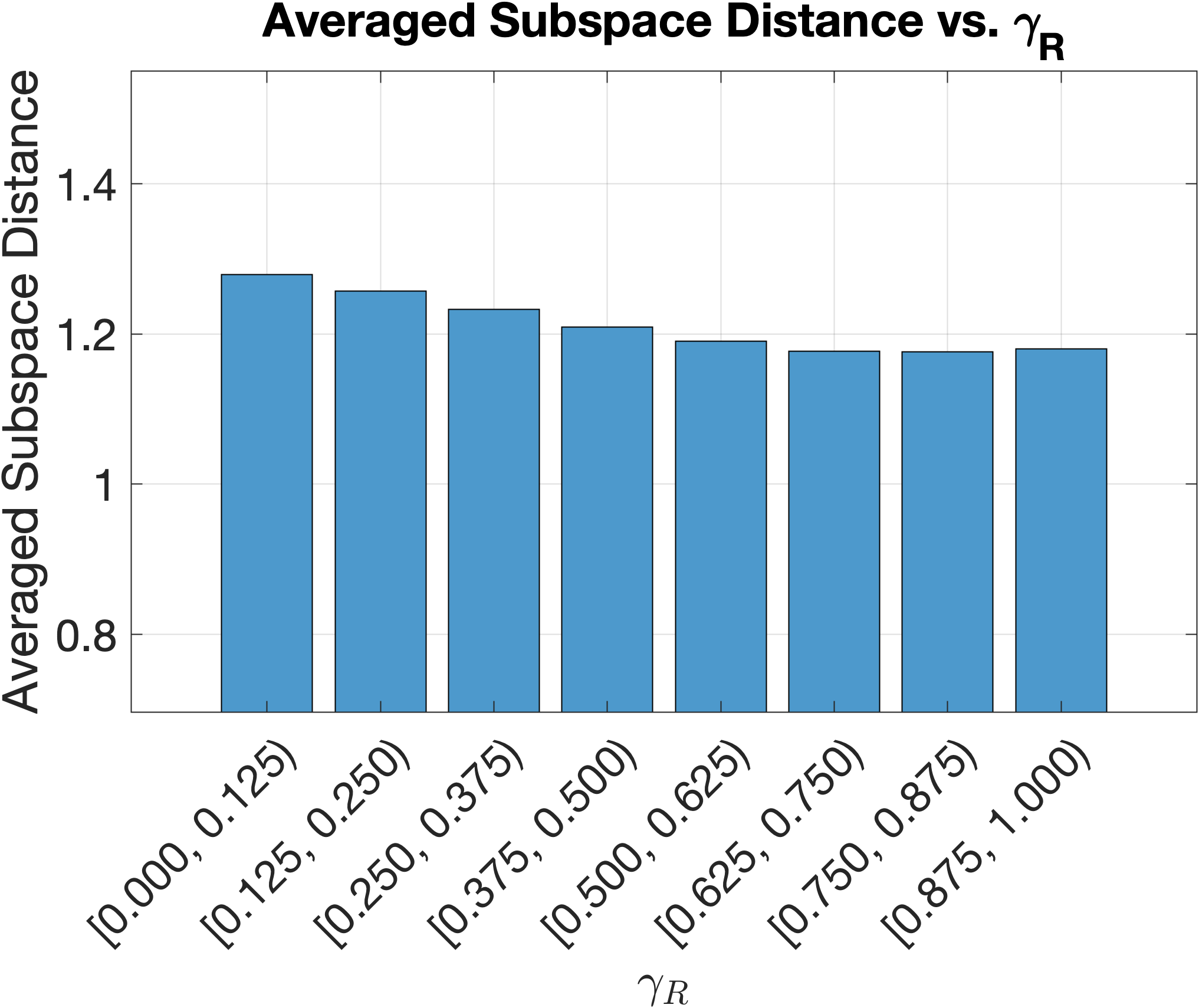}
    \caption{Averaged subspace distances by parameter in the Lotka-Volterra model.}
    \label{LV_sub_diff}
\end{figure}

The subspace distance is evaluated for each local region by analyzing the deviation in the local active subspace from the global subspace. To understand how these distances change as we move around the space, we report the average subspace distance for each parameter-subregion combination in Figure \ref{LV_sub_diff}.

From Figure \ref{LV_sub_diff}, we observe that the average subspace distances for $r_S$ and $r_R$ increase steadily as their values grow, with the lowest average distances occurring when each parameter is near zero. This indicates that the global and local active subspaces are most closely aligned when we are in a region of the parameter space where one of $r_S$ or $r_R$ is small, but that the distance between them is significantly larger when the growth rates are moderately large. Thus, global sensitivity rankings are more reliable in low-growth scenarios. Caution is warranted when applying global results to regions where the growth rates $r_S$ and $r_R$ are more robust; in such cases, the global rankings may not accurately reflect local parameter influence. 

The subspace distances for $\gamma_S$ and $\gamma_R$ show a milder decreasing trend as their values increase, suggesting that global sensitivity results are slightly less reliable in low-interaction regions when compared to high-interaction regions. In contrast, the subspace distances for $K_S$ and $K_R$ remain relatively stable and significant across their entire range. As such, there is no region that can be identified by the values of $K_S$ and $K_R$ alone in which we can be confident that global information accurately reflects local behavior. 

\begin{figure}[!htb]
    \centering    \includegraphics[width=.9\textwidth]{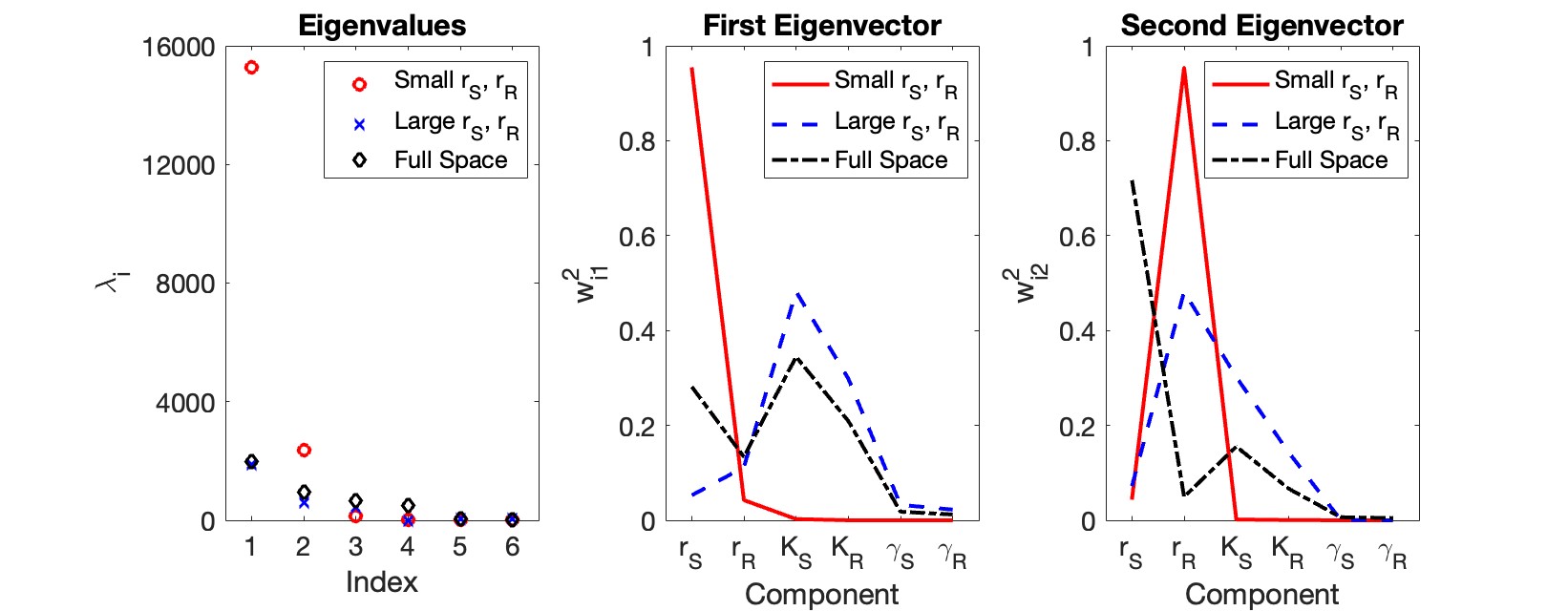}
    \caption{Eigendecomposition comparison for small growth, large growth, and full space scenarios.}
    \label{eigendecomp}
\end{figure}

While all three of the activity score, Morris elementary effects, and Sobol' indices methods reveal large discrepancies in sensitivity results across the full space, only the subspace distance plots used in the active subspace method enable us to determine where the issue is most prominent and refocus our analysis on a region in which sensitivity is more stable (i.e., a region in which we can trust that the globally-computed sensitivity metric is reflective of local behavior). In particular, the close alignment between global and local results in the case of small growth rates ($r_S, r_R\leq 0.125$) prompts us to examine this region further. Such analysis then uncovers why the globally-computed activity scores place such a focus on $r_S$ in addition to the two carrying capacities. This is illustrated in Figure \ref{eigendecomp}, where the eigenvalues and squared-eigenvector components are showcased for each of three scenarios: small growth only ($r_S, \ r_R \leq 0.25$ with all other parameters varying across their full range), large growth only ($r_S, \ r_R > 0.125$ with all other parameters varying across their full range), and the full admissible parameter space as used above. In the small growth case, the weight tips heavily in favor of $r_S$; this parameter is assigned the largest component in the first eigenvector and is paired with an overwhelmingly dominant eigenvalue. The second eigenvector prioritizes $r_R$, which---when paired with the second eigenvalue---then gets ranked as the second most influential parameter when the growth rates are restricted to $[0,0.125]$. The prioritization of $r_S$ over $r_R$ is due to the choice of initial $S$-to-$R$ ratio, $9:1$. 

In contrast, restricting the growth rates to $[0.125,1]$ results in a more even distribution of parameter influence, with the two carrying capacities $K_S$ and $K_R$ contributing the most and the growth rates playing a smaller role.  The full space analysis can be seen as an amalgamation of these two scenarios; combining behavior from both regions results in both $r_S$ and $K_S$ appearing highly influential in the global sense, even though the bulk of the influence for $r_S$ is derived only from the small growth rate scenario. 

Overall, the discrepancy between the global and local metrics reveals an important pitfall of relying solely on global sensitivity rankings for model reduction. For instance, a global analysis of this model using the active subspace method on the specified parameter space suggests that attention should be focused primarily on the growth rate and carrying capacity of Type-$S$ cells in subsequent tasks, while local analyses reveal that in many subregions, the Type-$S$ growth rate may actually be one of the least influential parameters. Sections \ref{sec:calibration}-\ref{sec:surrogate} will illustrate how the use of global sensitivity measures in the presence of such discrepancies may lead to issues with downstream tasks that rely upon model dimension reduction: specifically, model calibration and surrogate model construction.

\subsubsection{Model Calibration Illustration} \label{sec:calibration}

For models in which the full parameter set is not uniquely identifiable given the available data, it is common practice to conduct a sensitivity analysis and fix any uninfluential parameters at their nominal values prior to fitting the model to data. Here, we illustrate how using global information to select which parameters are fixed may lead to a less effective model calibration result. This is an example of a \textit{subset}-based reduction task, for which we are applying the activity score metric derived from the active subspace methodology.

\begin{figure}[!htb]
    \centering
    \includegraphics[width=\textwidth]{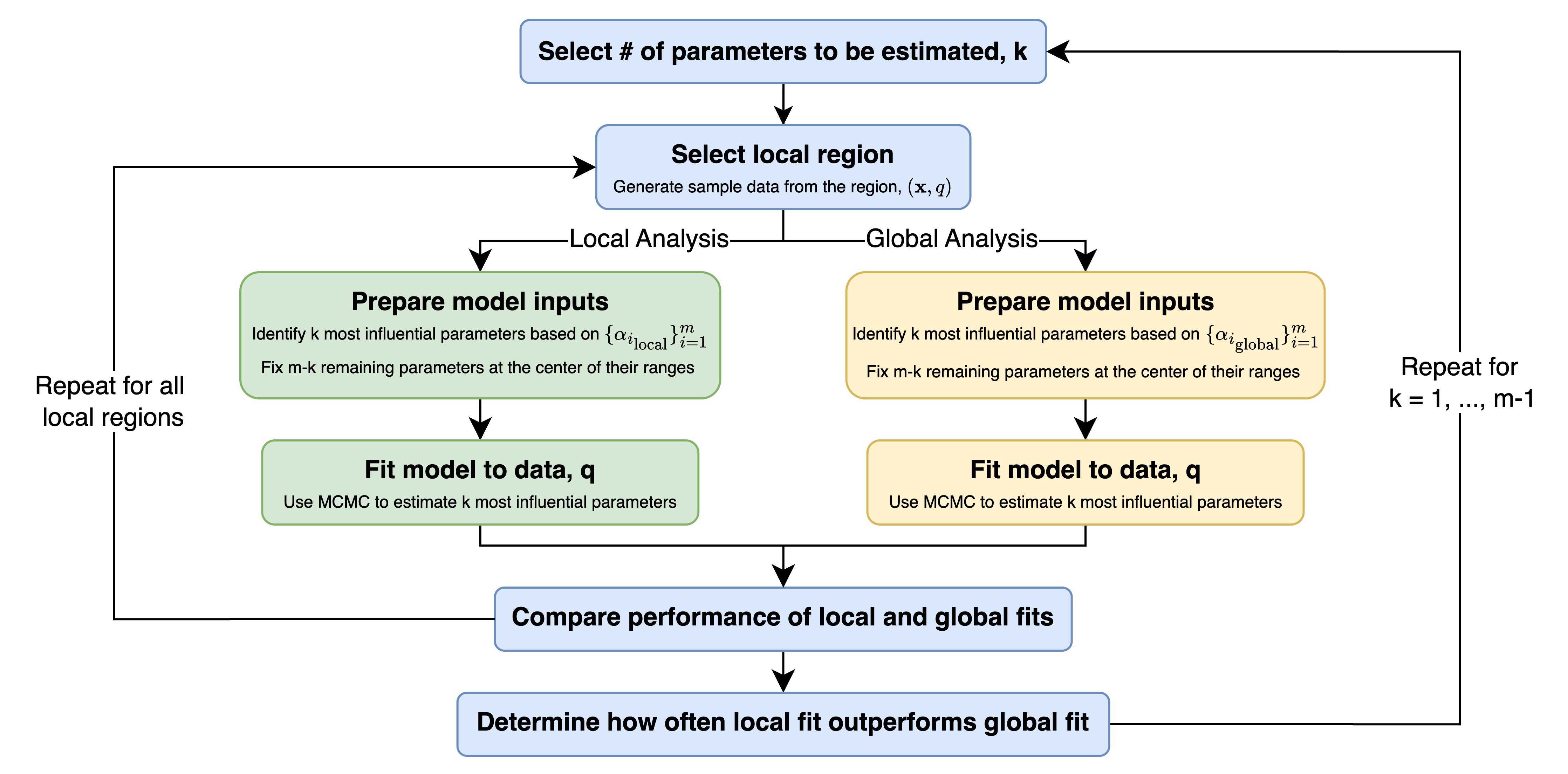}
    \caption{Model calibration experiment procedure.}
    \label{fig:calibrationdiagram}
\end{figure}

For each of the $8^6$ local regions defined previously for the Lotka-Volterra model, we generate a set of synthetic data using a parameter set drawn at random from that region. We then use the Delayed Rejection Adaptive Metropolis (DRAM) algorithm \cite{Haario, Smith}---a Markov chain Monte Carlo (MCMC) method for model calibration---to fit the model to our synthetic data in two ways: once, estimating only the top $k$ parameters identified by the global activity scores, and once estimating only the top $k$ parameters from the corresponding local activity score analysis, with the remaining $m-k$ parameters fixed at the center of their specified ranges. The performance of the resulting model fits (assessed using $(\text{QoI}_{\text{data}} - \text{QoI}_{\text{model}})^2$ as the optimized metric) is compared. This procedure is outlined in Figure \ref{fig:calibrationdiagram}. 

\begin{figure}[!htb]
\centering
\includegraphics[width=.5\textwidth]{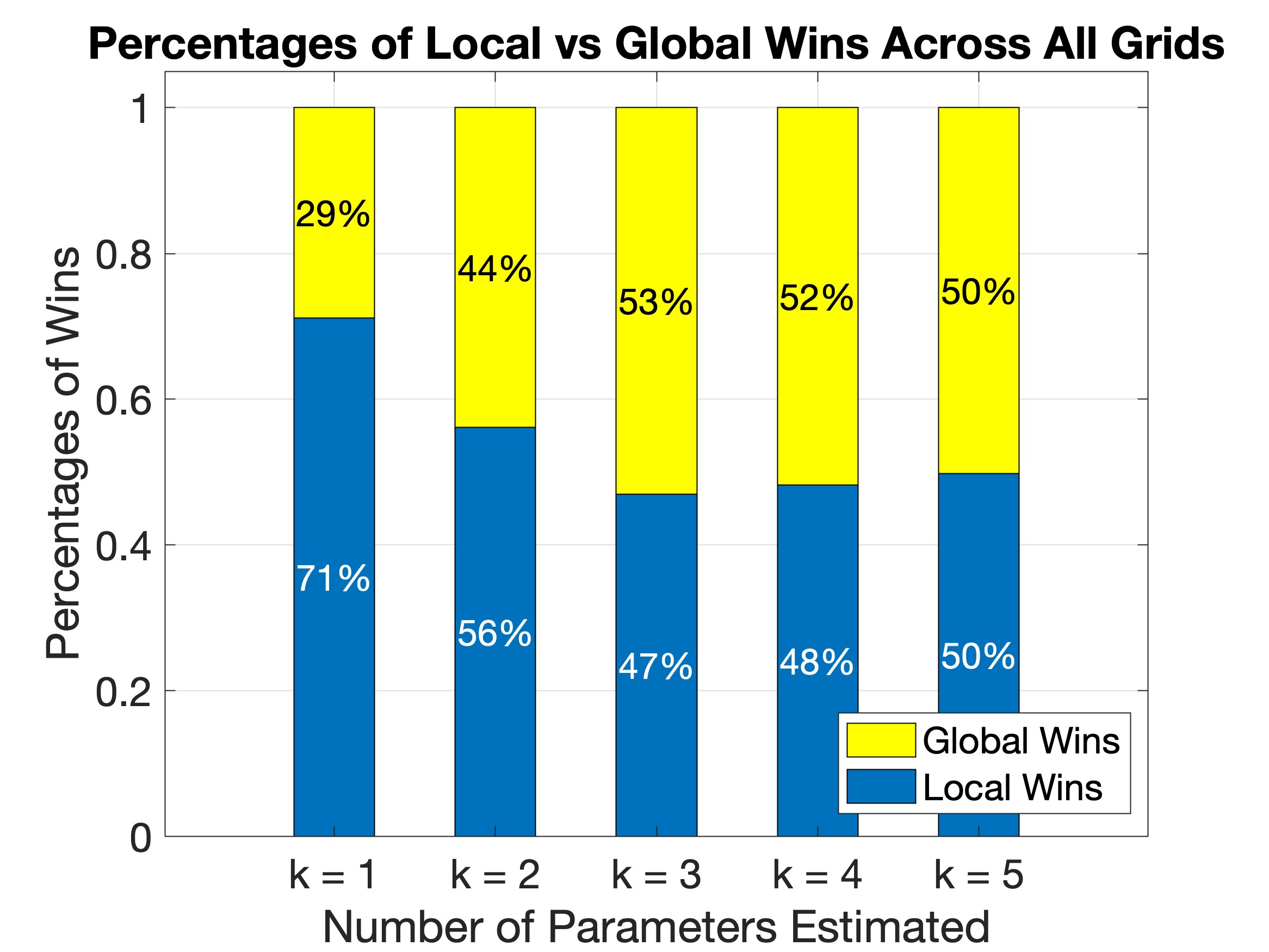}
\caption{Comparison of local versus global fit wins by dimension. }
\label{fig:LV_total_stacked}
\end{figure}

\begin{figure}[!htb]
\centering
\includegraphics[width=\textwidth]{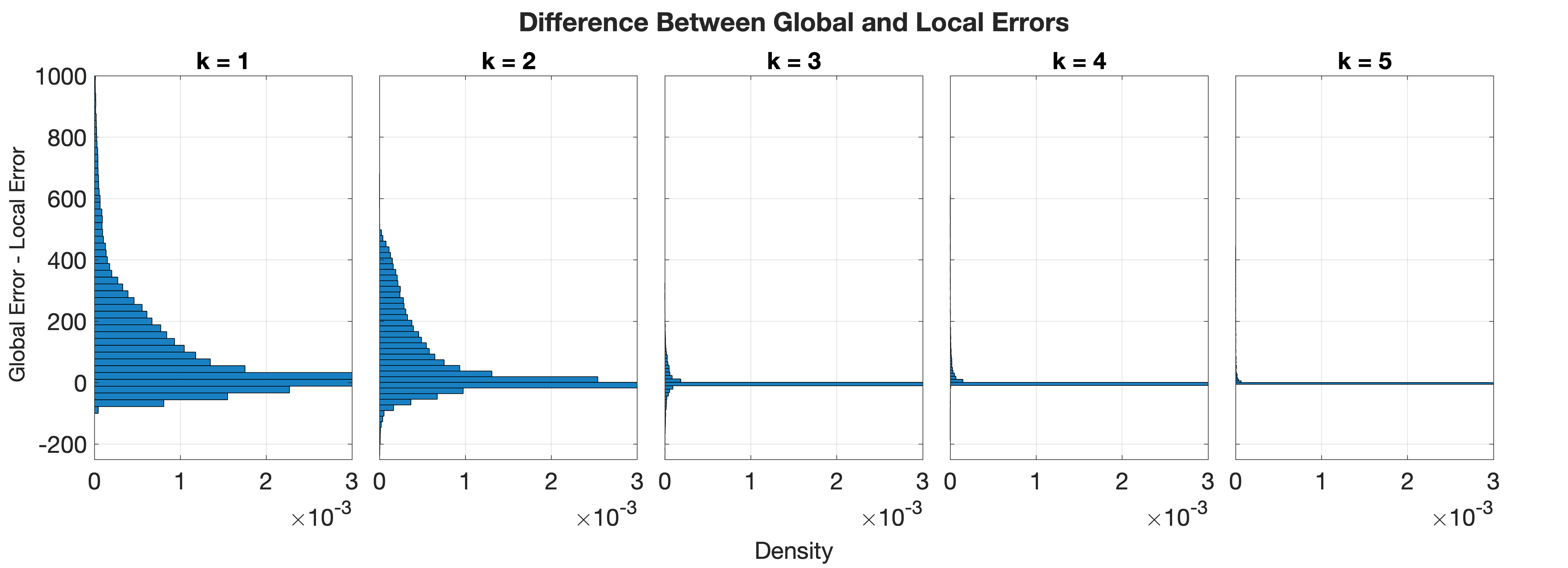}
\caption{Comparison of model fit errors resulting from global and local approaches by dimension. }
\label{fig:LV_errorcomp}
\end{figure}

Figures \ref{fig:LV_total_stacked}-\ref{fig:LV_errorcomp} show the results of our model calibration experiment. In Figure \ref{fig:LV_total_stacked}, we record the percentage of scenarios for which the local calibration results in the better model fit when compared to calibration using the global parameter ranking order. Specifically, calibration using local rankings dominates for $k = 1$ and $k = 2$; when only 1-2 parameters are estimated, it is crucial to focus attention on those that render the model most flexible in the region from which the data is produced to achieve the best possible fit, which requires the use of locally-sensitive parameters. As the number of estimated parameters increases, the advantage of the local approach narrows. With $k=3$, $k=4$, and $k=5$ parameters, the chosen subsets are often in agreement (or the additional parameters are insensitive and/or unidentifiable); thus, the model fits are nearly identical.  Figure \ref{fig:LV_errorcomp} illustrates the difference between the global fit errors and local fit errors. From this figure, we observe that in scenarios where the local calibration outperforms the global calibration, it tends to do so by a wider margin than in scenarios for which the global calibration edges out the local approach. This figure also confirms that in the $k=3$, $k=4$, and $k=5$ scenarios, the local and global model fits are nearly identical.

This experiment illustrates how the use of global sensitivity metrics may lead researchers to focus their attention on the wrong subset of parameters, reducing the efficacy of the calibration procedure. We suggest that if researchers observe large discrepancies between local and global active subspaces around the parameter space, they conduct further research to narrow the plausible ranges of their parameters prior to performing model calibration.

\subsubsection{Surrogate Model Construction Illustration} \label{sec:surrogate}

Another important task that relies heavily on sensitivity information is the construction of surrogate models to approximate the behavior of complex models, which may be too computationally expensive to evaluate repeatedly. Though the Lotka-Volterra model is not particularly expensive to evaluate, it serves as a useful illustrative example here. This experiment complements the earlier model calibration study by focusing on a \textit{subspace}-based reduction application, in which the directions to be disregarded are linear combinations of the original model parameters as defined by the inactive subspace. In particular, such information cannot be derived from our comparison methods (Morris elementary effects and Sobol' indices), but is readily available using active subspace methodology.

\begin{figure}[!htb]
    \centering
    \includegraphics[width=\textwidth]{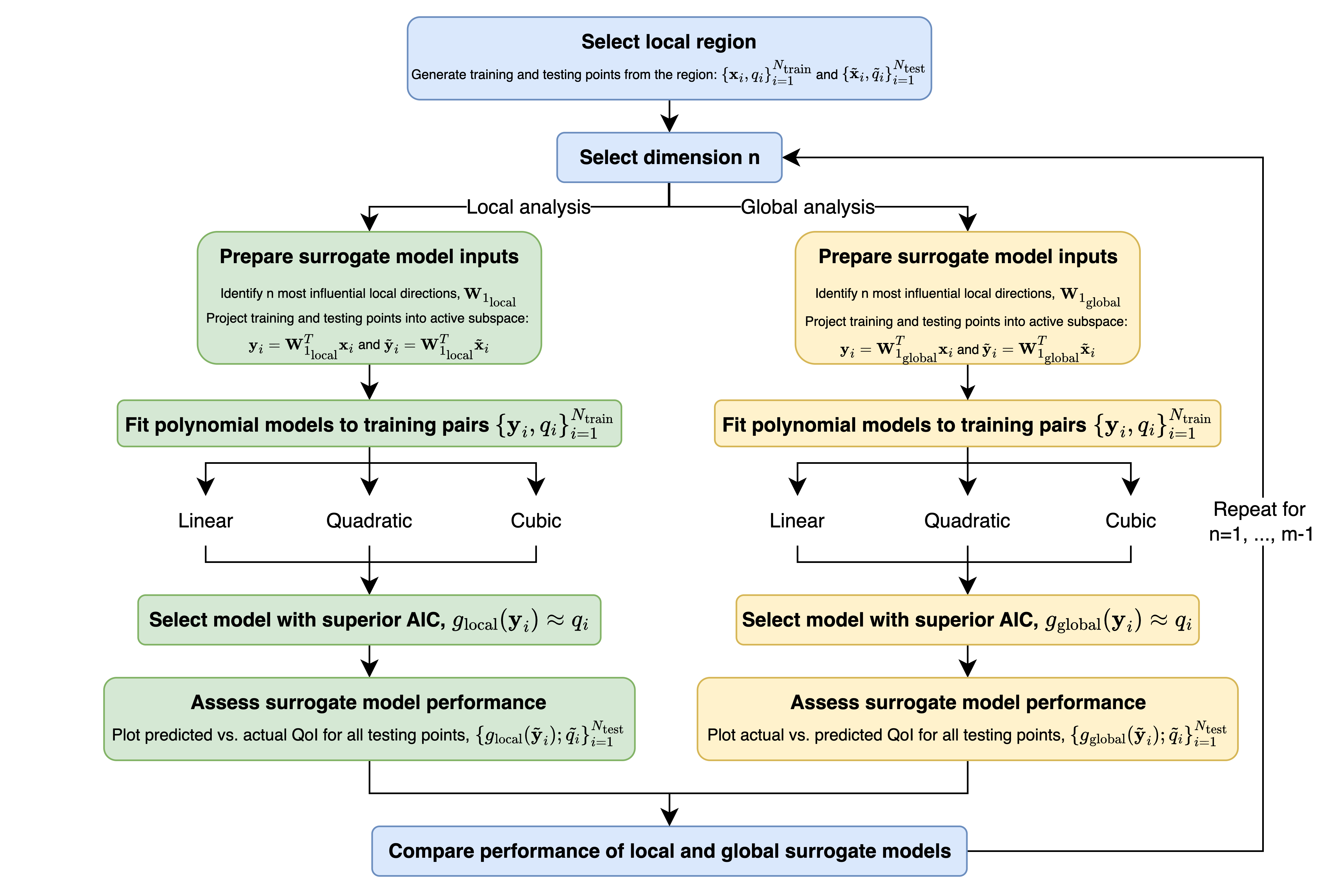}
    \caption{Surrogate modeling experiment procedure.}
    \label{fig:surrogatediagram}
\end{figure}

For this investigation, we construct polynomial response surface surrogate models on a subspace of $n$ dimensions, where $n<m$ (with $m$ representing the total number of input parameters). For a given local region, we build two surrogate models to approximate the QoI on that region: one using the globally calculated active subspace and one using the local active subspace. In each case, polynomial response surfaces of orders 1-3 are trained on points generated from the local region in question and projected into the active subspace. The ``best'' response surface model is chosen using the Akaike Information Criterion (AIC) \cite{Konishi}, where model performance is determined using a set of testing points generated independently from those used to train the models. Response surface performance is then compared between the two resulting models, with each being evaluated for its ability to reproduce the testing point outputs. The procedure is outlined in Figure \ref{fig:surrogatediagram}, and discussed in further detail in \cite{Col, Lewis, Smith}. 

\begin{figure}[!htb]
    \centering
    \includegraphics[width=\textwidth]{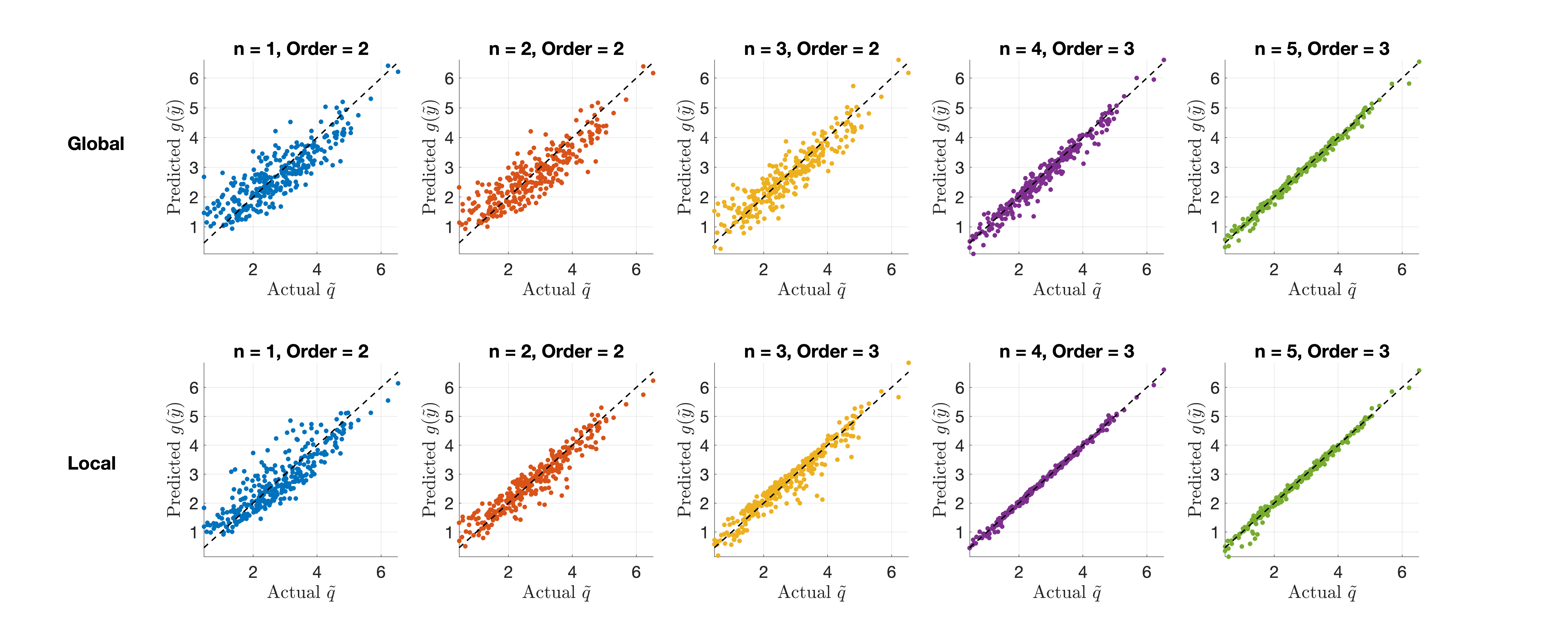}
    \caption{Predictive capability of global (top) versus local (bottom) response surfaces on the region $[0, 0.125]^6$.}
    \label{fig:sur_0}
\end{figure}

\begin{figure}[!htb]
    \centering
    \includegraphics[width=\textwidth]{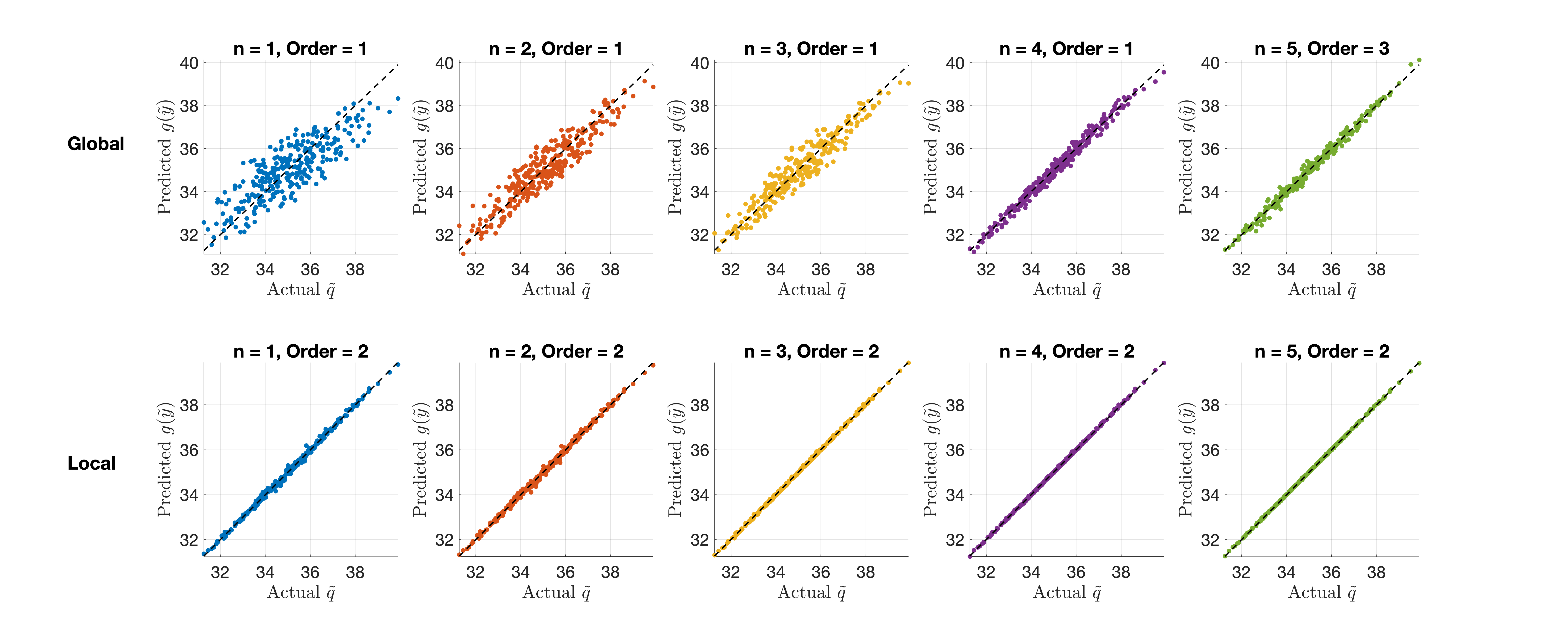}
    \caption{Predictive capability of global (top) versus local (bottom) response surfaces on the region $[0.5, 0.625]^6$.}
    \label{fig:sur_0.5}
\end{figure}

These comparisons are shown in Figures \ref{fig:sur_0} and \ref{fig:sur_0.5} for two of the $8^6$ local regions; note that other local regions display similar findings. On the local region where all six parameters are generated from the smallest bin, $[0, 0.125]^6$ (see Figure \ref{fig:sur_0}), we observe that the model constructed using local sensitivity information consistently outperforms the global model of the same input dimension, even on this subregion in which global results were determined to be in closest alignment with local results. Even more telling, the local model achieves a similar level of accuracy using two dimensions as the global model achieves with three. In this scenario, the three-dimensional global model---for which AIC selected a quadratic polynomial response surface---requires a total of 10 parameters to be estimated (with a model form of $c_1+c_2y_1+c_3y_2+c_4y_3+c_5y_1y_2+c_6y_1y_3+c_7y_2y_3+c_8y_1^2+c_9y_2^2+c_{10}y_3^2$, for active variables $y_1-y_3$ and parameter values $c_1-c_{10}$), in contrast to the two-dimensional local quadratic model, which requires only six (with a model form of $c_1+c_2y_1+c_3y_2+c_4y_1y_2+c_5y_1^2+c_6y_2^2$).  

In Figure \ref{fig:sur_0.5}, which shows the results for the more centralized region $[0.5,0.625]^6$---a region in which the discrepancy between local and global sensitivity information was more significant---the benefit of using local sensitivity information is even more apparent. The local model achieves a higher level of accuracy in a single dimension than the global model is able to achieve with five dimensions. Moreover, the local model built on this single dimension (with a quadratic response surface) requires the estimation of only three parameters (parameters $c_1-c_3$ for the model form $c_1+c_2y_1+c_3y_1^2$), while the global model built using a cubic polynomial on $n=5$ active variables requires the estimation of 56 parameters. The latter clearly no longer constitutes a simplification of the original model as would be desired.

In closing, we find that surrogate models built using global information with the goal of approximating local behavior yield less accuracy and tend to require larger dimensions to perform well. While this is generally not surprising---after all, local behavior will always be better approximated using local information---the decrease in global model accuracy (or increase in dimensions required to perform comparably to the local model) is compounded as the alignment in sensitivity information between local and global scales worsens. Thus, for models in which sensitivity varies widely across the parameter space, it is perhaps more beneficial to construct and piece together multiple surrogate models that are tuned to the local information from the regions on which they are constructed, rather than constructing a single global model that will suffer from reduced accuracy throughout the space.

\section{Conclusion} \label{sec:conclusion}

Global sensitivity metrics are commonly used in high-dimensional models to assess parameter importance. However, these metrics may misrepresent sensitivity when local variability is significant, leading to consequences in subsequent model analysis tasks. In this study, we utilize activity scores derived from the active subspace to demonstrate how parameter rankings may shift throughout the space, and quantify the divergence of local behavior from global results using the distance between locally- and globally-computed active subspaces. 

Using the Lotka-Volterra model, we find that large discrepancies exist between global and local parameter rankings across multiple established sensitivity methodologies, indicating that global summaries cannot adequately capture the full behavior of the system. These inconsistencies have numerous practical implications. A model calibration experiment shows that fitting the model using a subset of influential parameters chosen according to globally-derived sensitivity information often results in a poorer model fit, as compared to using local sensitivity information to select which parameters are to be estimated. A second experiment in surrogate model construction reveals that surrogate models that utilize global subspace information generally require higher-dimensional representations while yielding lower accuracy when compared to surrogate models constructed on locally-defined active subspaces. Such issues can result in flawed conclusions about parameter importance and may reduce the subsequent reliability of models in sensitive applications such as biological systems or engineering design.

Global sensitivity analysis is often heralded as a means by which to address a multitude of challenges related to model sensitivity and identifiability; however, we caution researchers to investigate the level of variability in local sensitivities for their model before using global sensitivity results to make consequential decisions. Active subspace methodology is well-suited for this purpose, since the active subspace can identify not just which parameters are most important in any given region, but also identify regions in which local and global sensitivity patterns are misaligned. Such analysis can then prompt suggestions for refinement of the parameter space to reduce variability prior to completing downstream tasks. 

Although conducting a thorough local analysis to accompany global measures may be computationally expensive, we believe that it is a worthwhile endeavor to avoid misrepresentation of parameter importance in subsequent analysis. In scenarios where a comprehensive local analysis is computationally prohibitive, we advise conducting this analysis on a random selection of local subregions as the computational budget allows, to assess whether sensitivity variability is likely to be an issue for the model under investigation. The development of strategies that adaptively switch between global and local metrics based on the stability of the active subspace---with the goal of reducing the computational expense of the proposed analysis---will be a focus of future investigations.

\end{document}